\documentclass[preprint]{imsart}
\usepackage{ifthen} % if else command
\usepackage{datetime} % Date and time
\settimeformat{ampmtime} % set in am-pm mode
\usepackage[usenames,dvipsnames,svgnames,table]{xcolor} % Using colours(including in  tables)
\usepackage{lastpage} % Access the last page of document (page x of lastpage)
\usepackage[final]{pdfpages}% Inserting external pdf pages.\includepdf command
\usepackage{graphicx,amsmath,amssymb}% importing graphics, maths, spl.symbols
\usepackage{multicol} % Multiple columns in the same document switch bw (single + double)
\usepackage{listings}% code listings, like verbatim package
\usepackage{bibunits}% For bibiliography, allows separate bibiliography for each unit
\usepackage{subfloat}% Including subcaption, subnumbering of figures, Fig1a,1b e.t.c
\usepackage{caption}
\usepackage{subcaption}%
\usepackage{stfloats}% Improves footnote position, skips e.t.c
\usepackage{eqlist}% listings with equation identations
\usepackage{eqparbox}% Required for making eqlist work, defining group of boxes with same width
\usepackage{pgf}% Making drawings and ploting in Latex(vector graphics)
\usepackage{xkeyval}%
\usepackage[colorlinks,citecolor=blue,urlcolor=blue]{hyperref}% typesetting of hyperlinks
\hypersetup{backref=true,colorlinks=true,raiselinks=true,hyperindex=true}
\usepackage{tikz} % Making drawings and ploting in Latex(vector graphics)
\usetikzlibrary{patterns,plothandlers,plotmarks,shapes,arrows,shapes.misc,backgrounds}
\usepackage{makeidx}% Creation of indexes (keywords and appearance)
\usepackage[title,titletoc]{appendix}% Making Appendix, subappendices
\usepackage[acronym,toc,style=list,nonumberlist]{glossaries}% Meanings of Abbreviations, op2- nomenclature package
\usepackage{fancyhdr}% change the default pagelayout
\usepackage{pdflscape}% Adds pdf supports, page layouts, works on driver level
\usepackage{cite}% Bibiliography and citations
\usepackage{xr}% Cross-referencing between different files
\usepackage{rotating}% Rotating any object with an arbitrary angle
\usepackage[below]{placeins}% Keeps floats in their place, preventing from moving past floatbarrier
\usepackage{url}% Adding urls
\usepackage{enumerate}% Different styles to enumerate environment
\usepackage{mathrsfs}% Different scripts for math symbols
\usepackage{amsthm}% Theorems, Algos and proofs
\usepackage{float}% placement of floats (figures,tables) at a place
\usepackage{multirow}% Tables with subrows and subcolumns
\usepackage[OT1]{fontenc}% Different encoding
\usepackage{exscale}% Big fonts especially maths
\usepackage[mathscr]{eucal}% Math calligraphy
\usepackage{stmaryrd}%Math calligraphy
\usepackage{verbatim}% like listings, but gives typewriter like impression
\usepackage{algorithm} % For algorithms
\usepackage{algpseudocode} % For algorithms
\usepackage[margin=1in]{geometry}% Handling Page Layout
\usepackage{lscape}
\usepackage{psfrag}
\usepackage{relsize}
\startlocaldefs
\newcommand{\SUM}[4]{%
{\displaystyle \sum_{{#1}={#2}}^{#3}\!{\!#4}
}%
}
\newcommand{\SUML}[2]{{\displaystyle \sum_{#1}\!{\!#2}
}}
\newcommand{\PROD}[4]{%
{\displaystyle \prod_{{#1}={#2}}^{#3}\!{\!#4}
}%
}
\newcommand{\PRODP}[3]{{\displaystyle \prod_{#1}^{#2}\!{\!#3}
}}
\newcommand{\NrmOne}[1]{\|{#1}\|_1}
\newcommand{\NrmTwo}[1]{\|{#1}\|_2}
\newcommand{\NrmZ}[1]{\|{#1}\|_0}
\newcommand{\INT}[3]{%
{\displaystyle \int\limits_{#1}^{#2}\!{\!#3}
}%
}
\newcommand{\argmin}[1]{\underset{#1}{\operatorname{arg}\,\operatorname{min}}\;} % use as \argmin{x}
\numberwithin{equation}{section}
\theoremstyle{plain}
\newtheorem{Thm}{Theorem}[section]
\newtheorem{Def}{Definition}[section]
\newtheorem{Cor}{Corollary}[section]
\newtheorem{Exp}{Example}[section]
\newenvironment{EqA}[2]
{\begin{equation}\label{#1}\begin{array}{#2}}
{\end{array}\end{equation}}
\newenvironment{DefNoNm}[1]
{\begin{Def}\label{#1}}
{\par\end{Def}}
\newenvironment{ThmNoNm}[1]
{\begin{Thm}\label{#1}}
{\par\end{Thm}}
\newenvironment{CorNoNm}[1]
{\begin{Cor}\label{#1}}
{\par\end{Cor}}
\newenvironment{ExpNoNm}[1]
{\begin{Exp}\label{#1}}
{\par\end{Exp}}
\newcommand{\SeE}{\end{subequations}}
\newcommand{\SeS}{\begin{subequations}}
\newcommand{\iu}{{i\mkern1mu}}
\def\b{\mathbf{b}}
\def\c{\mathbf{c}}

\def\e{\mathbf{e}}

\def\g{\mathbf{g}}
\def\h{\mathbf{h}}
\def\gm{\boldsymbol{\gamma}}
\def\m{\mathbf{m}}

\def\r{\mathbf{r}}
\def\s{\mathbf{s}} % Unknown Source Parameter in b = As
\def\u{\mathbf{u}}
\def\v{\mathbf{v}}
\def\w{\mathbf{w}}
\def\x{\mathbf{x}}
\def\y{\mathbf{y}}
\def\xhat{\hat{\mathbf{x}}}
\def\xhats{\hat{x}}
\def\btil{\tilde{\mathbf{b}}}
\def\br{\breve{\mathbf{b}}}

\def\zhat{\hat{\mathbf{z}}}
\def\zhats{\hat{z}}
\def\A{\mathbf{A}}
\def\B{\mathbf{B}}
\def\D{\mathbf{D}}

\def\H{\mathbf{H}}
\def\I{\mathbf{I}}

\def\R{\mathbf{R}}

\def\W{\mathbf{W}}
\def\mean{\mathbb{E}}

\def\pr{\mathbb{P}}
\def\re{\mathbb{R}}
\def\sgn{\mathbb{S}}

\def\cf{\mathbb{F}}
\def\rl{\mathbb{R}}

\def\ocf{\mathbb{G}}
\def\cI{\mathscr{I}}
\def\cJ{\mathscr{J}}
\def\pJ{\mathscr{P}}
\def\eJ{\mathscr{E}}
\endlocaldefs
\begin{document}
%---------------------------------------------
\begin{frontmatter}
%-------------------------------------------------------------------------
\title{The LASSO Estimator: Distributional Properties}
%\thanksref{T1}}
\runtitle{The LASSO Estimator: Distributional Properties}
%\thankstext{T1}{Footnote to the title with the ``thankstext'' command.}
\begin{aug}
\author{\fnms{Rakshith} \snm{Jagannath}\ead[label=e1]{rejarr@gmail.com}}
\and
\author{\fnms{Neelesh S} \snm{Upadhye}\ead[label=e2]{neelesh@iitm.ac.in}}
\runauthor{Jagannath and Upadhye}

\affiliation{Indian Institute of Technology, Madras}

\address{Indian Institute of Technology, Madras\\
Sardar Patel Road, Adyar, Chennai,\\
Tamil Nadu-600036\\
\printead{e1}\\
\phantom{E-mail:\ }\printead*{e2}}
\end{aug}
%--------------------------
\begin{abstract}
The least absolute shrinkage and selection operator (LASSO) is a popular technique for simultaneous estimation and model selection. There have been a lot of studies on the large sample asymptotic distributional properties of the LASSO estimator, but it is also well-known that the asymptotic results can give a wrong picture of the LASSO estimator's actual finite-sample behavior. The finite sample distribution of the LASSO estimator has been previously studied for the special case of orthogonal models. The aim in this work is to generalize the finite sample distribution properties of LASSO estimator for a real and linear measurement model in Gaussian noise.

In this work, we derive an expression for the finite sample characteristic function of the LASSO estimator, we then use the Fourier slice theorem to obtain an approximate expression for the marginal probability density functions of the one-dimensional components of a linear transformation of the LASSO estimator.
\end{abstract}
%-------------------------
%--------------------------------------------
\end{frontmatter}
%---------------------------------------------
%-------------------------------------------------------
\section{Introduction and Motivation}\label{intro}
LASSO (least absolute shrinkage and selection operator) \cite{tibs-lasso,chen-lasso} has been developed as a tool to find sparse solutions of the linear regression problem. It has been used extensively in an expanding field of applications from statistics to estimation scenarios with remarkably good results. It is used extensively in the parameter estimation framework to estimate the unknown parameter(s) with guarantees for the estimation error (model fit) \cite{CS_Baraniuk,donoho-cs}. As a robust approximation of the well known Maximum Likelihood (ML) estimator, the LASSO can be realized robustly and efficiently through convex Second Order Cone (SOC) programming techniques \cite{Ben-Tal-LMC} and unlike the efficient subspace methods \cite{Arraysp}, the LASSO technique is reliable even with one data measurement realization (single snapshot) \cite{complexlars}.

The regularization parameter (sparsity threshold parameter) in the LASSO is a mathematical tool to implement the compromise between model fit and the estimated model order, which is the number of nonzero entries in the estimate. As the regularization parameter evolves, the LASSO solution changes continuously, forming a continuous trajectory in a very high dimensional space which is referred to as the LASSO path \cite{Efronlars,complexlars}.

The LASSO algorithm, in general assumes the knowledge of sparsity of the unknown parameter or the optimum regularization parameter for estimating the unknown parameter. However, these are generally unknown and need to be estimated.

In the detection framework, the focus is to propose detection tests to estimate the optimal sparsity threshold parameter, so that the number of non-zero entries (or the sparsity) and their corresponding locations (indices) in the estimate is same as the actual parameter. The performance evaluation of the detection tests is done using the $p$-values or the probability of correct detection. As sparsity plays an important role in the estimation performance, this issue has been recognized as a significant gap between theory and practice by several authors \cite{cvcs, sure-eldar}. The problem of estimating the sparsity of the parameter from the measurement data is fundamental to many other applications such as estimating the critical number of measurements for successful recovery, design of the sensing matrix e.t.c \cite{lopes-CS}. The exploration of the detection framework of the LASSO is fairly recent (e.g. \cite{siglass} and its citations). Sparsity and model order estimation techniques like statistical cross-validation, Mallow’s Cp selection, Stein's unbiased risk estimator and Bayesian information criteria (BIC) and its variants have been proposed in \cite{cs_cv,Mosesequivalence,lopes-CS} for estimation of sparsity (or sparsity threshold parameter). Bayesian based LASSO estimation, wherein the linear model is interpreted from a Bayesian perspective and the sparsity threshold parameter is modeled as a hyper-prior for estimation, have been proposed in \cite{FL}. Hence we see that an accurate estimate of the sparsity or the sparsity threshold parameter is critical for enhancing the performance of parameter estimation using LASSO. However, the performance of the above techniques depends on the initial guess of the sparsity threshold parameter or the choice of a grid for the sparsity threshold parameter. These techniques also use the distribution of the ML estimator instead of the distribution of the LASSO estimator for deriving the tests for model order estimation.
Hence, we note that the detection framework requires the distribution of the LASSO estimator for proposing reliable detection tests and evaluating their performance. So, there is a need for studying the distributional properties of the LASSO estimator. Hence, in this work we explore the distribution of the LASSO estimator.

There has been a lot of related work in the understanding of the asymptotic distributional properties of LASSO estimator, e.g., \cite{knight2000}, but it is also well-known that the asymptotic results can give a wrong picture of the LASSO estimator's actual finite-sample behavior \cite{paul-critic,potscher-critic,potscher-finite}. In particular, \cite{potscher-finite} studies the finite sample LASSO distribution for the case of orthogonal models and shows that the asymptotic results do not provide a reliable assessment for the finite sample distribution. Therefore there is a need for studying the finite sample distributional properties of the LASSO estimator. Hence, in this work,  we will study the finite sample characteristic function (cf) of the LASSO estimator for any general model matrix.

This article is organized as follows. In Section-\ref{prelim}, we discuss the notations and state without proof, some well known theorems that will be used in this work. In Section-\ref{results}, we give the details of our results on the cf and approximate probability density function (pdf) of the LASSO estimator. In Section-\ref{sim}, we perform numerical simulations to verify the results discussed in Section-\ref{results}. In Section-\ref{concl}, we discuss the conclusions and some possibilities for future work. We give the details of the proofs of the theorems stated in Section-\ref{results} in  Appendix-\ref{proofs}.
%-------------------------------------------------------
\section{Preliminaries}\label{prelim}
In this section, we discuss the notations, measurement model, LASSO estimator and some well known theorems used in this work.
\subsection{Notations}\label{notations}
We use bold lower case letters to represent vectors ($\x$), bold upper case letters to represent matrices ($\A$) and scripted letters to represent sets ($\cJ$, generally finite index sets). $\cJ\setminus\cI$ denotes the set-difference operation between sets $\cJ$ and $\cI$. For a given matrix (vector) $\A$, $\A^T$ denotes the regular transpose, $|\A|$ and $\A^{-1}$ denote the determinant and inverse of the square matrix $\A$ and $\A^{\dagger}$ denotes the Moore-Penrose pseudo inverse of $\A$. For a vector $\x$, $\NrmZ{\x}$, $\NrmOne{\x}$, $\NrmTwo{\x}$ denote the $l_0$ pseudo norm which is equal to the number of non-zero elements in $\x$, $l_1$ and $l_2$ norms respectively. $\mean$ denotes expectation, $\pr$ denotes probability. For a random vector $\x$, $f_{\x}(\x)$ denotes its pdf. Convolution of $f([x_{1},x_{2},\ldots,x_{N}])$ with $g(x_{1})$ is denoted equivalently by $f(\x)\star g(x_{1})$ or $f(\overline{x}_1,x_2,\ldots,x_N)$ or $f(\overline{x}_{1},\x^{-})$, which is defined as the following integral
%-------------------------------------------------
%------------------------------------
\begin{EqA}{One/p1/Eqn/Ea1}{rcl}
%-----------------------------------
f(\overline{x}_{1},\x^{-})=f(\x)\star g(x_{1})&=& \INT{\re}{}{f(x_{1}-u,\x^{-})g(u)du}
%-----------------------------------
\end{EqA}
%-----------------------------------
%-------------------------------------------------
where $\x^{-} = [x_2,\dots,x_N]$. Similarly $f(\overline{\x}_{\cI},\x_{\cJ\setminus\cI})$ denotes the $|\cI|$ dimensional extension of convolution across the dimensions of $\x$ given by $\cI$ and finally $\sgn(\mathbf{u})$ denotes the sign of elements of $\mathbf{u}$, we have
%-------------------------------------------------
%-----------------------------------
\begin{EqA}{One/p1/Eqn/Ea2}{rcl}
%-----------------------------------
\sgn(u_i) & = & \left\{
%-----------------------------------
\begin{array}{rl}
-1 & \text{if } u_i < 0,\\
1 & \text{if } u_i > 0.
%-----------------------------------
\end{array} \right.
%-----------------------------------
\end{EqA}
%-----------------------------------
$\sgn(u_i)$ is arbitrary for $u_i = 0$.

Let $f:\rl^{N}\rightarrow \rl$ be a function. We define the following operations on $f$.
\begin{enumerate}
\item \textbf{Integral Projection:} $\mathcal{I}_{\x}$ is the projection operator that reduces an $N$ dimensional function to $M$ dimensions by integrating out $N-M$ dimensions:
    \begin{equation*}
    \mathcal{I}_{\x}[f](x_1,x_2,\ldots,x_M) = \INT{}{}{f(x_1,\ldots,x_N) dx_{M+1}\ldots dx_N}
    \end{equation*}
\item \textbf{Slicing:} $\mathcal{S}_{\x}$ is the slicing operator that reduces an $N$ dimensional function to an $M$ dimensional function, by zeroing out $N-M$ dimensions: $\mathcal{S}_{\x}[f](x_1,\ldots,x_M) = f(x_1,\ldots,x_M,0,\ldots,0)$.
\item \textbf{Change of Basis:} Let $\B$ denote a full rank $N\times N$ matrix and let $\x$ denote an $N$ dimensional vector. Then $\B[f](\x) = f(\B^{-1}\x)$.
\end{enumerate}
The action of multiple operators on the function is denoted by $\odot$. For example, $\mathcal{I}_{\B\x}\odot\B^{-1}[f]$ denotes the change of basis followed by projection operation on the function $f(\x)$.

\subsection{Measurement Model}
We consider the following linear regression model,
%-------------------------------------------------
\begin{EqA}{One/p1/Eqn/Ea3}{rcl}
\b &=& \A\x + \v,
\end{EqA}
%-----------------------------------
where $\b = [b_1,b_2,\ldots,b_M]$ denotes the measurement vector of length $M$, $\A$ denotes the model matrix of size $M\times N$, $\v$ denotes the white Gaussian noise with zero mean and covariance matrix, $\sigma^2\I$ and $\x = [x_1,x_2,\ldots,x_N]$ denotes the sparse parameter vector of length $N$ and sparsity $K$ (number of non-zero entries in $\x$), which needs to be estimated. We also assume that columns of $\A$ have unit norm. We focus on the LASSO estimator \cite{tibs-lasso,chen-lasso}, which estimates the sparse parameter $\x$ in (\ref{One/p1/Eqn/Ea3}) by solving the following convex optimization problem,
%-------------------------------------------------
%------------------------------------
\begin{EqA}{One/p1/Eqn/Ea4}{rcl}
%-----------------------------------
\xhat &=& \argmin{\x} \tau\NrmOne{\x} + \frac{1}{2}\NrmTwo{\b-\A\x}^2.
%-----------------------------------
\end{EqA}
%-----------------------------------
where $\tau>0$ is the sparsity thresholding parameter, which controls the sparsity of $\xhat$. We observe that (\ref{One/p1/Eqn/Ea4}) is a convex relaxation of the following combinatorial problem
%-------------------------------------------------
%------------------------------------
\begin{EqA}{One/p1/Eqn/Ea5}{rcl}
%-----------------------------------
\min\NrmZ{\x}&\text{subject to}&\hspace{1mm} \NrmTwo{\b-\A\x}\leq\sigma
%-----------------------------------
\end{EqA}
%-----------------------------------
It has been shown in \cite{candes-rip} that the LASSO estimator is an exact relaxation of (\ref{One/p1/Eqn/Ea5}) if the model matrix, $\A$ satisfies the restricted isometry property (R.I.P) defined below in Definition-\ref{One/p1/Defns/Def1}, and hence gives the sparsest estimate to the linear regression problem of (\ref{One/p1/Eqn/Ea3}) (see Theorem-\ref{One/p1/Thm/Thm1} below).
%-------------------------------------------------
%------------------------
\begin{DefNoNm}{One/p1/Defns/Def1}
%-----------------------------------
For each integer $K = 1, 2,\ldots$, we define the restricted isometry constant $\delta_K$ of a matrix $\A$ as the smallest number such that
%-----------------------------------
%------------------------------------
\begin{EqA}{One/p1/Eqn/Ea6}{rcl}
%-----------------------------------
(1-\delta_K)\NrmTwo{\y}^2&\leq&\NrmTwo{\A\y}^2\leq(1+\delta_K)\NrmTwo{\y}^2
%-----------------------------------
\end{EqA}
%-----------------------------------
holds for all $K$ sparse vectors, $\y$. A vector is said to be $K$ sparse if it has at most $K$ nonzero entries.
%-----------------------------------
\end{DefNoNm}
%-----------------------------------
%------------------------
\begin{ThmNoNm}{One/p1/Thm/Thm1}
%-----------------------------------
Assume that $\delta_{2K}<\sqrt{2}-1$ and $\NrmTwo{\v}<\sigma$. Then the solution $\xhat$ to (\ref{One/p1/Eqn/Ea4}) obeys
%-----------------------------------
%------------------------------------
\begin{EqA}{One/p1/Eqn/Ea7}{rcl}
%-----------------------------------
\NrmTwo{\xhat-\x}&\leq&C_0 \frac{1}{\sqrt{K}}\NrmOne{\x-\x^{K}}+C_1\sigma
%-----------------------------------
\end{EqA}
%-----------------------------------
%-----------------------------------
for some constants $C_0$ and $C_1$. In particular, if $\x$ is $K$-sparse, the recovery is exact. Here $\x^{K}$ is the best sparse approximation one could obtain if one knew exactly the locations and amplitudes of the $K$-largest entries of $\x$.
\vspace{-2.0ex}
%-----------------------------------
\begin{proof}
See \cite{candes-rip}
\end{proof}
%-----------------------------------
\end{ThmNoNm}
%-----------------------------------
We now state some well known definitions and theorems required for this work.
%-------------------------------------------------
%------------------------
\begin{DefNoNm}{One/p1/Defns/Def2}
%-----------------------------------
 The cf of a random vector $\y$ is,
%-----------------------------------
%------------------------------------
\begin{EqA}{One/p1/Eqn/Ea8}{rcl}
%-----------------------------------
\cf_{\y}(\u) &=& \mean_{\y}\{\exp(\iu\u^{T}\y)\}\\
&=& \INT{-\infty}{\infty}{f_{\y}(\y)\exp(\iu\u^{T}\y)d\y}
%-----------------------------------
\end{EqA}
%-----------------------------------
%-----------------------------------
Clearly, cf is the Fourier transform of the pdf of $\y$ with $\u$ as the variable in the Fourier domain. The cf has the properties like, cf is a uniformly continuous, bounded and hermitian function with guaranteed existence, $\cf_{\y}(\mathbf{0}) = 1$ and cf is a bijection with probability distributions, i.e, for any two random variables $X_1$ and $X_2$, both have the same probability distribution if and only if $\cf_{X_1}= \cf_{X_2}$.
%-----------------------------------
\vspace{-0.2ex}
%-----------------------------------
\end{DefNoNm}
%-----------------------------------
%--------------------------------------------------------------
\begin{ThmNoNm}{One/p1/Thm/Thm2}
%--------------------------------------------------------------
A Borel probability measure $\pr$ on $\rl^{N}$ is uniquely determined by its one dimensional projections, i.e. a probability measure on Euclidean space is uniquely determined by the values it gives to half-spaces.
%-----------------------------------
\begin{proof}
See \cite{cwold}
\end{proof}
%-----------------------------------
\end{ThmNoNm}
%-----------------------------------
%------------------------
\begin{ThmNoNm}{One/p1/Thm/Thm3}
%-----------------------------------
Let $f$ be an $N$ dimensional function, let $\cf$, $\B$, $\mathcal{I}$,
$\mathcal{S}$ represent the Fourier transform, change of basis, projection and slicing operations as explained above, then we have
%-----------------------------------
%------------------------------------
\begin{EqA}{One/p1/Eqn/Ea9}{rcl}
%-----------------------------------
\cf\odot\mathcal{I}_{\B^{-1}\x}\odot\B[f] &=& \mathcal{S}_{\B^{T}\u}\odot\frac{\B^{-T}} {|\B^{-T}|}\odot\cf[f]
%-----------------------------------
\end{EqA}
%-----------------------------------
\vspace{-0.5cm}
%-----------------------------------
\begin{proof}
See \cite{Ng}
\end{proof}
%-----------------------------------
\end{ThmNoNm}
%-------------------------------------------------
Next, we discuss the main results of this work in the next section.
%-------------------------------------------------------
\section{Main Results}\label{results}
%-------------------------------------------------------
In this section, we first derive an expression for the cf of the LASSO estimator. The expression for cf is appropriately sliced to derive the one dimensional projections of a linear transformation of the LASSO estimator. The one dimensional projections yield the approximate marginal pdfs of components of the linear transformation of the LASSO estimator.
%-------------------------------------------------
%------------------------
\begin{ThmNoNm}{One/p1/Thm/Thm4}
%-----------------------------------
The cf of the LASSO estimator, $\cf_{\xhat}(\u)$ as a function of the true parameter, $\x$ is given by the following implicit relationship
%-----------------------------------
\begin{EqA}{One/p1/Eqn/Ea10}{rclll}
%-----------------------------------
\cf_{\mathsmaller{\W\xhat + \tau\sgn(\xhat)}}(\u)&=& \SUML{\cI\in\pJ}{\Big(\PRODP{\substack{k\in\cI\\j\in\cJ\setminus\cI}}{}{\sin(\tau u_k)\cos(\tau u_j)}\Big)\cf_{\xhat}(\overline{\c}_{\cI},\c_{\cJ\setminus\cI})}&=&\exp\big(\iu\u^{T}\W\x-\frac{\sigma^2}{2}\u^{T}\W\u\big)
%-----------------------------------
\end{EqA}
%-----------------------------------
where $\cJ = \{1,2,\ldots,N\}$ is an index set, $\pJ$ is the power set of $\cJ$, $\cI$ be an element of $\pJ$ or equivalently $\cI\subseteq\cJ$, $\W = \A^{T}\A$, $\c_{\cI} = \u^{T}\w_j,j\in\cI$ similarly $\c_{\cJ\setminus\cI} = \u^{T}\w_k,k\in\cJ\setminus\cI$ and $\cf_{\xhat}(\overline{\c}_{\cI},\c_{\cJ\setminus\cI})$ denotes the $|\cI|$ dimensional convolution of $\cf_{\xhat}(\c_{\cI},\c_{\cJ\setminus\cI})$ with $g(\c_{\cI}) = \PRODP{j\in\cI}{}{\hspace{1mm}\frac{-1}{\pi c_j}}$.
%-----------------------------------
\vspace{-2.0ex}
%-----------------------------------
\end{ThmNoNm}
%-----------------------------------
%-------------------------------------------------
Example-\ref{One/p1/Thm/Eg1} below illustrates Theorem-\ref{One/p1/Thm/Thm4}
%-------------------------------------------------
%------------------------------------------
\begin{ExpNoNm}{One/p1/Thm/Eg1}
%-----------------------------------------
When $\xhat$ has $N = 3$ elements, we have $\cJ = \{1,2,3\}$, the power set, $\pJ = \{\emptyset,{1},{2},{3},\{1,2\},\{2,3\},\{1,3\},\{1,2,3\}\}$ and $\cI$ is one of the elements of $\pJ$. Hence, the relationship between $\cf_{\xhat}(\u)$ and $\x$ is given by
%-----------------------------------
%-------------------------------------
\begin{EqA}{One/p1/Eqn/Ea11}{rcl}
%------------------------------------
\cf_{\mathsmaller{\W\xhat+\tau\sgn(\xhat)}}(\u)&=&\cos(\tau u_1)\cos(\tau u_2)\cos(\tau u_3)\cf_{\xhat}(c_1,c_2,c_3)+ \sin(\tau u_1)\cos(\tau u_2)\cos(\tau u_3)\cf_{\xhat}(\overline{c}_1,c_2,c_3)+\\
&&\cos(\tau u_1)\sin(\tau u_2)\cos(\tau u_3)\cf_{\xhat}(c_1,\overline{c}_2,c_3)+
\cos(\tau u_1)\cos(\tau u_2)\sin(\tau u_3)\cf_{\xhat}(c_1,c_2,\overline{c}_3)+\\
&&\sin(\tau u_1)\sin(\tau u_2)\cos(\tau u_3)\cf_{\xhat}(\overline{c}_1,\overline{c}_2,c_3)+
\cos(\tau u_1)\sin(\tau u_2)\sin(\tau u_3)\cf_{\xhat}(c_1,\overline{c}_2,\overline{c}_3)+\\
&&\sin(\tau u_1)\cos(\tau u_2)\sin(\tau u_3)\cf_{\xhat}(\overline{c}_1,c_2,\overline{c}_3)+
\sin(\tau u_1)\sin(\tau u_2)\sin(\tau u_3)\cf_{\xhat}(\overline{c}_1,\overline{c}_2,\overline{c}_3)\\
&=&\exp\big(\iu\u^{T}\W\x-\frac{\sigma^2}{2}\u^{T}\W\u\big)\nonumber
%-----------------------------------
\end{EqA}
%-----------------------------------
%-----------------------------------
\vspace{-2.0ex}
%-----------------------------------
\end{ExpNoNm}
%-----------------------------------
where $\cf_{\xhat}(\overline{c}_{j},\ldots,)$ represents convolution as defined in (\ref{One/p1/Eqn/Ea1}) with $g(c_j) = \frac{-1}{\pi c_j}$
%-------------------------------------------------

\textbf{Remarks:} We make the following observations from Theorem-\ref{One/p1/Thm/Thm4}.
\begin{enumerate}
\item We observe that the relationship between $\cf_{\mathsmaller{\W\xhat+\tau\sgn(\xhat)}}(\u)$ and $\cf_{\xhat}(\c = \W\u)$ has $2^N$ terms in total, which does not simplify the evaluation of the pdf of LASSO estimator except for the special case of orthogonal model matrix ($\W= \I$), where the entries of $\xhat$ become independent and hence the pdf of $\xhat$ can be obtained easily (see Corollary-\ref{One/p1/Thm/Thm5}). Proof of Corollary-\ref{One/p1/Thm/Thm5} (Appendix-\ref{prf2}) is an alternate way (\cite{potscher-finite}) to obtain pdf of $\xhat$, when $\W$ is orthogonal.
\item We observe that the one dimensional projections of the cf can be evaluated from (\ref{One/p1/Eqn/Ea10}) by slicing.
\item We then invoke the Cramer-Wold theorem (Theorem-\ref{One/p1/Thm/Thm2}) to conclude that the one-dimensional projections are sufficient for the evaluation of the joint distribution of the LASSO estimator.
\item We note that for $\tau = 0$, we obtain the cf relationship for the maximum likelihood (ML) estimator as a special case and hence the corresponding pdf of the ML estimator \cite{Kayesti} can be easily evaluated.
\end{enumerate}

Next, we evaluate the one dimensional projections of the cf, $\cf_{\xhat}$ for different cases of $\W$. We first start with the case of $\W = \I$.
%-------------------------------------------------
%------------------------
\begin{CorNoNm}{One/p1/Thm/Thm5}
%-----------------------------------
If $\W = \I$, where $\I$ is the identity matrix, or any diagonal matrix $\D$ (in general), then the estimates $\hat{x}_i, i=1,2,\ldots,N$ are all independent and hence the cf of the lasso estimator is given by,
%-----------------------------------
%------------------------------------
\begin{EqA}{One/p1/Eqn/Ea12}{rcl}
%-----------------------------------
\cf_{\xhats_{k}}(u_{k})\cos(\tau u_{k})+ \cf_{\xhats_{k}}(\overline{u}_{k})\sin(\tau u_{k})&=& \exp\big(\iu u_{k}x_{k}-\frac{\sigma^2}{2}u^{2}_{k}\big), k = 1,2,\ldots,N
%-----------------------------------
\end{EqA}
%-----------------------------------
%-----------------------------------
and the marginal pdf of the individual components, $\hat{x}_i, i=1,2,\ldots,N$ is obtained by simply applying the inversion theorem to (\ref{One/p1/Eqn/Ea12}) as,
%-----------------------------------
%------------------------------------
\begin{EqA}{One/p1/Eqn/Ea13}{rcl}
%-----------------------------------
f_{\xhats_k}(\xhats_k)& = & \left\{
%-----------------------------------
\begin{array}{rl}
\frac{1}{\sqrt{2\pi\sigma^2}}\exp\big(-\frac{(\xhats_k+\tau-x_k)^2}{2\sigma^2}\big)& \text{if } \xhats_k>0,\\
\frac{1}{\sqrt{2\pi\sigma^2}}\exp\big(-\frac{(\xhats_k-\tau-x_k)^2}{2\sigma^2}\big) & \text{if } \xhats_k<0.
%-----------------------------------
\end{array} \right.
%-----------------------------------
\end{EqA}
%-----------------------------------
%-----------------------------------
\end{CorNoNm}
%-----------------------------------
%-------------------------------------------------
\textbf{Remarks:}
\begin{itemize}
\item The expression (\ref{One/p1/Eqn/Ea13}) is similar to the expression ($5$) of \cite{potscher-finite} (we note that \cite{potscher-finite} have derived the pdf of the estimation error, $\xhats-x$).
\item We note that the joint pdf of $\xhat$ is just the product of the marginal pdfs of $\xhats_{k}, k=1,2,\ldots,N$ obtained in (\ref{One/p1/Eqn/Ea13}).
\end{itemize}

Now, we evaluate the cf for the case when $\W$ is a full rank matrix.
%-------------------------------------------------
%------------------------
\begin{CorNoNm}{One/p1/Thm/Thm6}
%-----------------------------------
Let $\W$ be any full rank matrix, $\zhat = \W\xhat$, $\w_j$ be the $j^{th}$ column of $\W$ and $\zhats_j$ be the $j^{th}$ element of $\zhat$. Then the cf of the one dimensional projections, $\zhats_k$ for components $k$ corresponding to large $|\xhats_{k}|$ can be approximated as,
%-----------------------------------
%------------------------------------
\begin{EqA}{One/p1/Eqn/Ea16}{rcl}
%-----------------------------------
\cos(\tau u)\cf_{\zhats_k}(u)+ \sin(\tau u)\cf_{\zhats_k}(\overline{u})
&=& \exp(-u^{2}\frac{\sigma^2}{2}w_{kk})\exp\Big(ju(\w_{k}^{T}\x)\Big)
%-----------------------------------
\end{EqA}
%-----------------------------------
%-----------------------------------
and the marginal pdfs of the components $\hat{z}_k$, corresponding to the large $|\xhats_{k}|$ is obtained by simply applying the inversion theorem to (\ref{One/p1/Eqn/Ea16}) as,
%-----------------------------------
\begin{EqA}{One/p1/Eqn/Ea17}{rcl}
%-----------------------------------
f_{\zhats_k}(\zhats_k)& = & \left\{
%-----------------------------------
\begin{array}{rl}
\frac{1}{\sqrt{2\pi\sigma^2w_{kk}}}\exp\big(-\frac{(\zhats_k+\tau-\w^{T}_{k}\x)^2}{2\sigma^2}\big)& \text{if } \zhats_k>0,\\
\frac{1}{\sqrt{2\pi\sigma^2w_{kk}}}\exp\big(-\frac{(\zhats_k-\tau-\w^{T}_{k}\x)^2}{2\sigma^2}\big) & \text{if } \zhats_k<0.
%-----------------------------------
\end{array} \right.
%-----------------------------------
\end{EqA}
%-----------------------------------
\end{CorNoNm}
%-------------------------------------------------
\textbf{Remarks:}
\begin{enumerate}
\item In Corollary-\ref{One/p1/Thm/Thm6}, we make the approximation $\sgn(\h_{k}^{T}\zhat)\approx \sgn(\zhats_k)$ to evaluate the expression, $\mathcal{S}_{\u}[\cf_{\xhat}(\W\u)\star(\frac{-1}{\pi c_k})]$ (more details in Appendix-\ref{prf3}). Here $\h_{k}$ is a column of $\H = \W^{-1}$. This approximation is valid whenever $|h_{kk}\zhats_{k}|>>\SUM{i}{1,i\neq k}{N}{|h_{ki}\zhats_i|}$. Since $\W$ is a symmetric positive definite (pd) matrix, $\H$ is also symmetric and pd, hence the diagonal entries of $\H$ are positive. So, if $\H$ is diagonally dominant and the index $k$ corresponds to a large entries of $|\xhats_{k}|$, then the approximation works well. From the simulations (Section-\ref{sim}) also, we observe that the approximation works well for the index $k$ corresponding to the large non-zero entries of $|\xhats_{k}|$.
\item The evaluation of exact expression for $\mathcal{S}_{\u}[\cf_{\xhat}(\W\u)\star(\frac{-1}{\pi c_k})]$ requires some prior knowledge or assumptions on the pdf $f_{\mathsmaller{\zhat}}(\zhat)$. In Appendix-\ref{hilbertval}, we derive an expression for $\mathcal{S}_{\u}[\cf_{\xhat}(\W\u)\star(\frac{-1}{\pi c_k})]$ by assuming the inherent distribution of $f_{\mathsmaller{\zhat}}(\zhat)$ to be multivariate Gaussian. This assumption is justified because the  expression in (\ref{One/p1/Eqn/Ea10}), which represents some operations on cf of $\xhat$ is equal to cf of Gaussian pdf. Also, for the special case of orthogonal $\W$, we obtain a Gaussian pdf for $f(\xhat)$. Hence, we can make the assumption that $f_{\mathsmaller{\zhat}}(\zhat)$ is multivariate Gaussian. From the resulting expression for $\mathcal{S}_{\u}[\cf_{\xhat}(\W\u)\star(\frac{-1}{\pi c_k})]$, we then justify the use of the approximation $\sgn(\h_{k}^{T}\zhat)\approx \sgn(\zhats_k)$ for $k$ corresponding to large $|\xhats_{k}|$.
\end{enumerate}

Further, we evaluate the cf for any general case of $\W$.
%-------------------------------------------------
%------------------------
\begin{CorNoNm}{One/p1/Thm/Thm7}
%-----------------------------------
Let $\W$ be any general matrix, $\zhat = \W\xhat$, $\w_j$ be the $j^{th}$ column of $\W$ and $\zhats_j$ be the $j^{th}$ element of $\zhat$. Then the cf of the one dimensional projections, $\zhats_k$ for components $k$ corresponding to large $|\xhats_{k}|$ can be approximated as,
%-----------------------------------
\begin{EqA}{One/p1/Eqn/Ea16c}{rcl}
%-----------------------------------
\cos(\tau u)\cf_{\zhats_k}(u)+ \sin(\tau u)\cf_{\zhats_k}(\overline{u})
&=& \exp(-u^{2}\frac{\sigma^2}{2}w_{kk})\exp\Big(ju(\w_{k}^{T}\x)\Big)
%-----------------------------------
\end{EqA}
%-----------------------------------
and the marginal pdfs of the components $\hat{z}_k$, corresponding to the large $|\xhats_{k}|$ is obtained by simply applying the inversion theorem to (\ref{One/p1/Eqn/Ea16c}) as,
%-----------------------------------
\begin{EqA}{One/p1/Eqn/Ea17}{rcl}
%-----------------------------------
f_{\zhats_k}(\zhats_k)& = & \left\{
%-----------------------------------
\begin{array}{rl}
\frac{1}{\sqrt{2\pi\sigma^2w_{kk}}}\exp\big(-\frac{(\zhats_k+\tau-\w^{T}_{k}\x)^2}{2\sigma^2}\big)& \text{if } \zhats_k>0,\\
\frac{1}{\sqrt{2\pi\sigma^2w_{kk}}}\exp\big(-\frac{(\zhats_k-\tau-\w^{T}_{k}\x)^2}{2\sigma^2}\big) & \text{if } \zhats_k<0.
%-----------------------------------
\end{array} \right.
%-----------------------------------
\end{EqA}
%-----------------------------------
\end{CorNoNm}
%-------------------------------------------------
\textbf{Remarks:}
\begin{enumerate}
\item In Corollary-\ref{One/p1/Thm/Thm7}, we again make the approximation $\sgn(\h_{k}^{T}\zhat)\approx \sgn(\zhats_k)$ to evaluate $\mathcal{S}_{\u}[\cf_{\xhat}(\W\u)\star(\frac{-1}{\pi c_k})]$ (more details in Appendix-\ref{prf3}). Now, $\h_{k}$ is a column of $\H = \W^{\dagger}$. The approximation is valid whenever $|h_{kk}\zhats_{k}|>>\SUM{i}{1,i\neq k}{N}{|h_{ki}\zhats_i|}$. Since $\W$ is a symmetric positive semi definite (psd) matrix, $\H$ is also symmetric and psd, hence the diagonal entries of $\H$ are non-negative. So, if the matrix is diagonally dominant and the index $k$ corresponds to the non-zero entry of $\xhat$, then the approximation works well. From the simulations also we observe that the approximation seems to work for the index $k$ corresponding to non-zero entries of $\xhat$.
\item Again the evaluation of exact expression for $\mathcal{S}_{\u}[\cf_{\xhat}(\W\u)\star(\frac{-1}{\pi c_k})]$ requires some prior knowledge or assumptions on the pdf $f_{\mathsmaller{\zhat}}(\zhat)$ and the exact expression confirms the feasibility of the approximation $\sgn(\h_{k}^{T}\zhat)\approx \sgn(\zhats_k)$ for the components $k$ corresponding to large $|\xhats_{k}|$.
\end{enumerate}
%-------------------------------------------------------
\section{Numerical Simulations}\label{sim}
%-------------------------------------------------------
In this section, we perform simulations to validate the theoretical results for cf. The simulations are performed using the linear regression model of (\ref{One/p1/Eqn/Ea3}). In the simulation setup, we choose a measurement vector, $\b$ of length $M = 4$. The model matrix $\A$ is chosen as a Hadamard matrix for orthogonal case and random matrix of size $M$ by $N$ consists of entries from Bernoulli distribution for non-orthogonal case. The columns of the model matrix $\A$ are always normalized to have unit norm. $N = 4$ for orthogonal and full rank models and $N = 8$ for singular model. The parameter vector $\x$ is chosen such that its sparsity $K=1$ and its $\frac{N}{2}+1$ entry is non-zero
($|x_{\mathsmaller{N/2+1}}| = 4$ when $N=4$ and $|x_{\mathsmaller{N/2+1}}| = 8$ when $N=8$). The noise is generated as a multivariate Gaussian random vector noise covariance variance matrix $\sigma^2\I$, where $\sigma = 1$. In the following, we use Monte-Carlo simulations for $L = 10000$ noisy realizations to obtain $L$ noisy estimate vectors. We then run the LASSO algorithm using the MATLAB CVX package \cite{cvx} to solve the optimization problem of (\ref{One/p1/Eqn/Ea4}) to obtain the LASSO estimate $\xhat$ for each realization. We choose $\tau = 1$ for orthogonal and full rank models and $\tau = 2$ for singular model.
%------------------------------
\begin{figure}[h]
%\vspace{-2ex}
%\FigRule
\begin{tabular}{llll}
%--------------------------
\hspace{-2ex}
\begin{subfigure}[t]{0.3\textwidth}
\includegraphics[width=4cm,height=5cm]{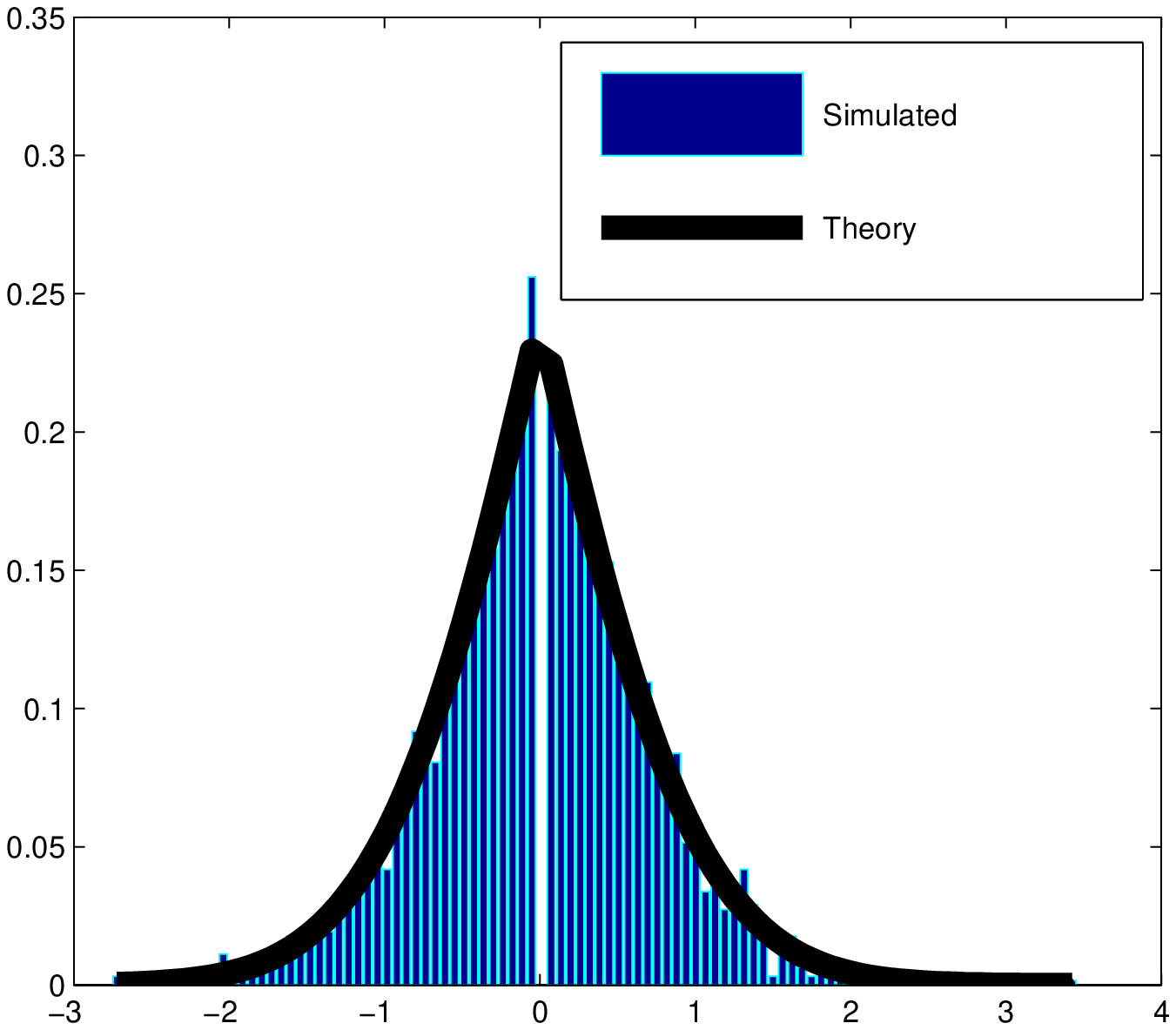}
\caption{$\hat{x}_1$}
\label{oneorder8ant}
\end{subfigure}
%--------------------------
&
%-----------------------------
\hspace{-10ex}
\begin{subfigure}[t]{0.3\textwidth}
\includegraphics[width=4cm,height=5cm]{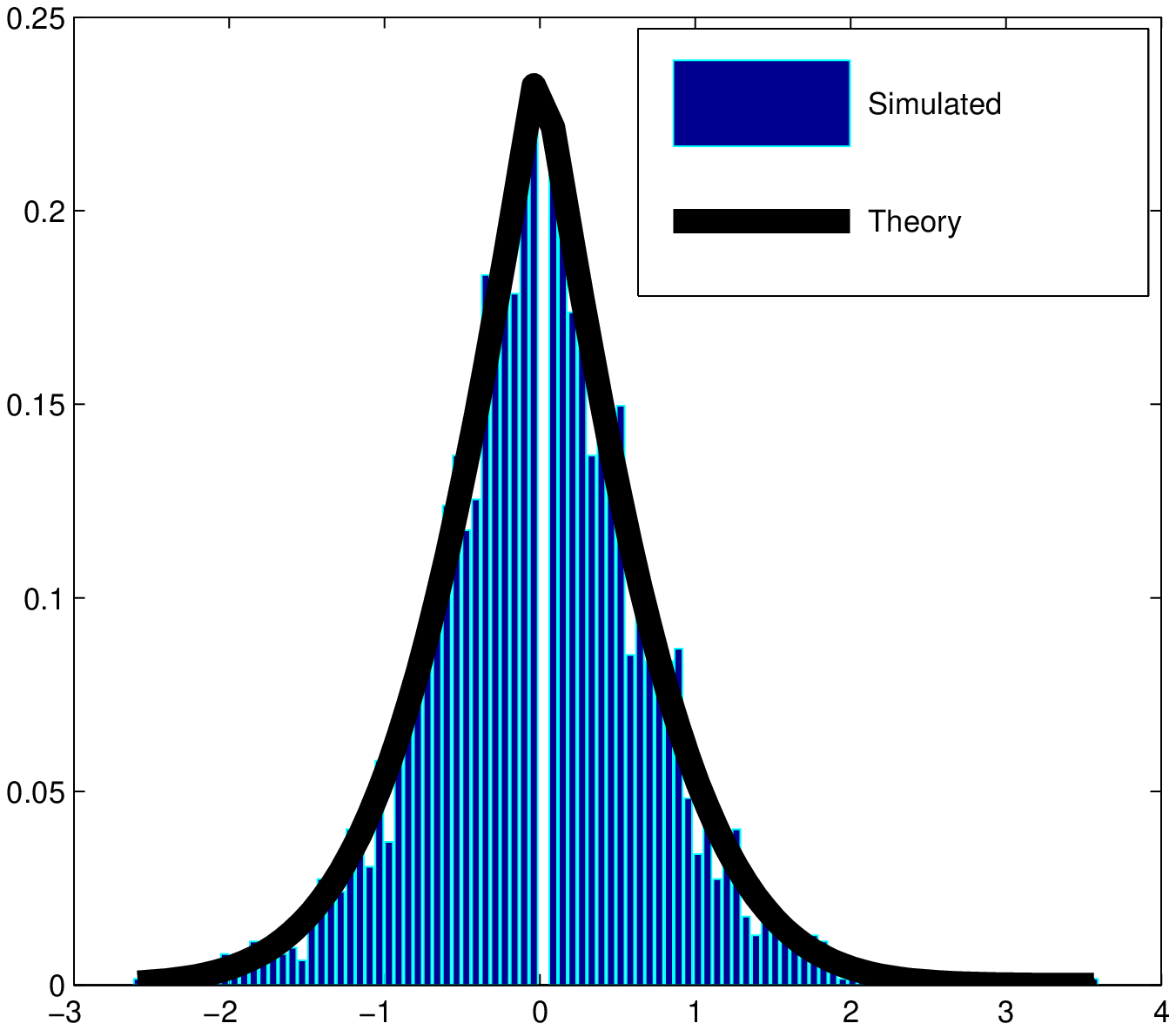}
\caption{$\hat{x}_2$}
\label{oneorder16ant}
\end{subfigure}
%-----------------------------
&
%-----------------------------
\hspace{-8ex}
\begin{subfigure}[t]{0.3\textwidth}
\includegraphics[width=4cm,height=5cm]{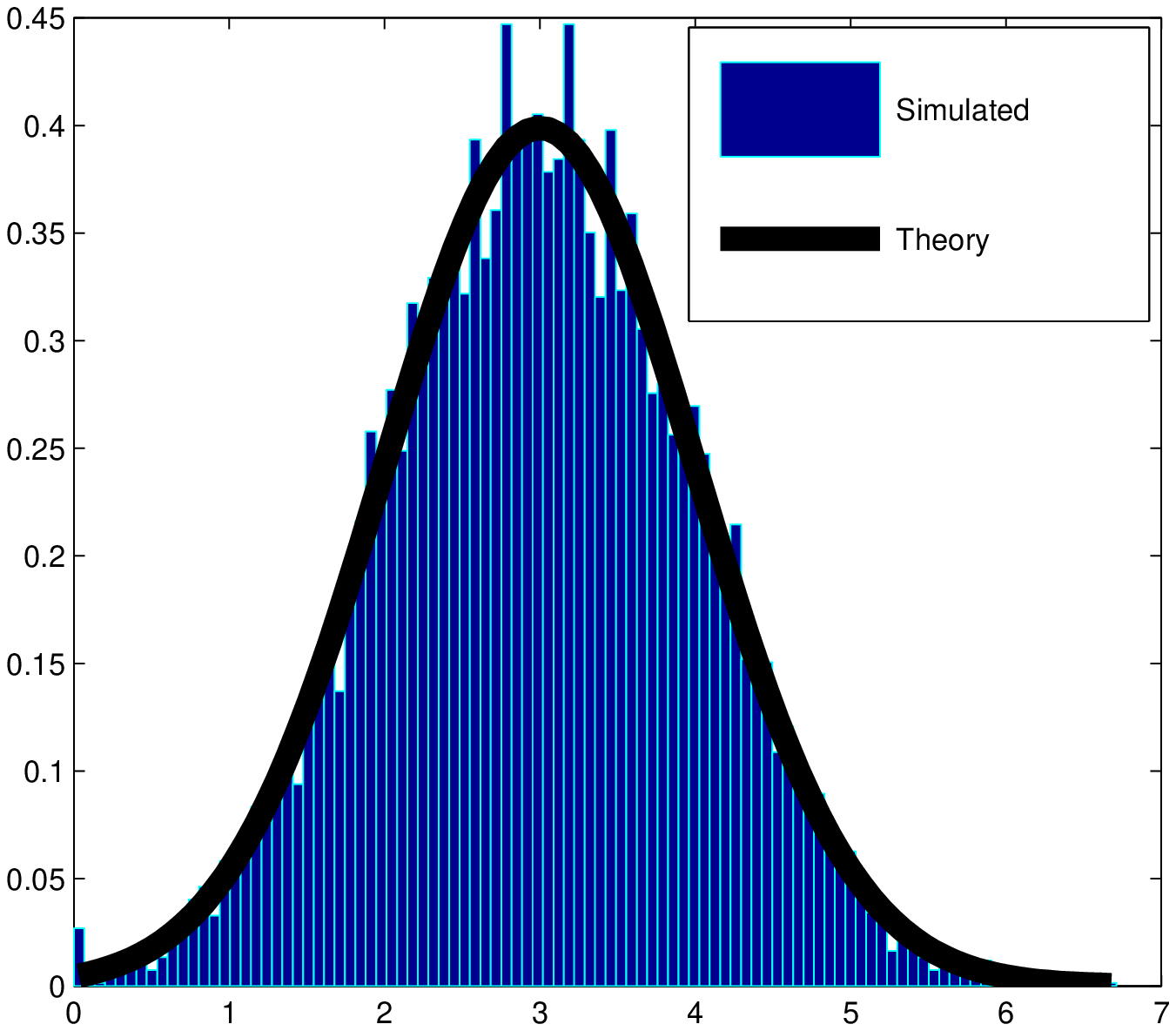}
\caption{$\hat{x}_3$}
\label{twoorder8ant}
\end{subfigure}
%-----------------------------
&
%-----------------------------
\hspace{-10ex}
\begin{subfigure}[t]{0.3\textwidth}
\includegraphics[width=4cm,height=5cm]{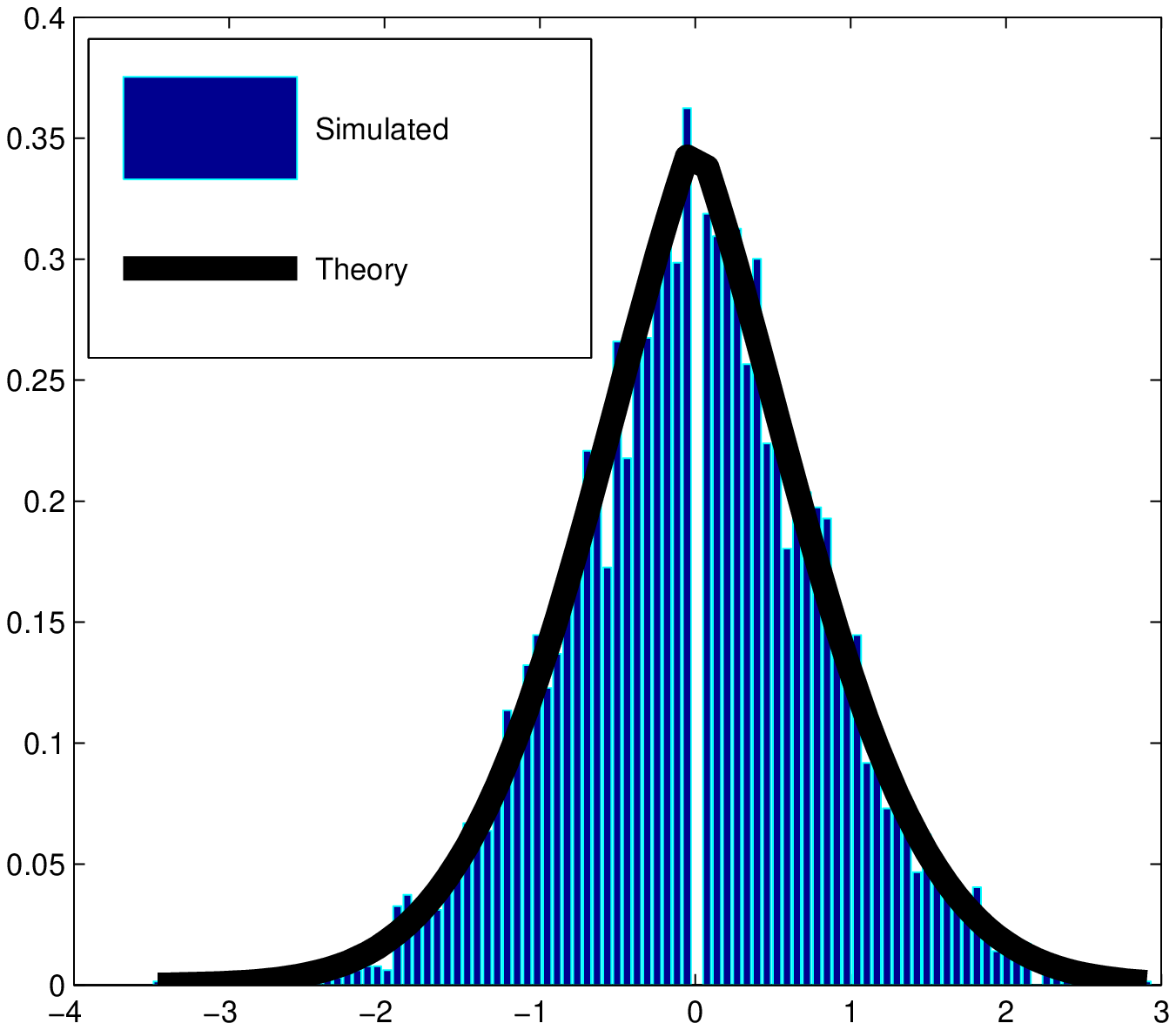}
\caption{$\hat{x}_4$}
\label{x4hat}
\end{subfigure}
%----------------------------
\end{tabular}
\caption{pdf of the LASSO estimator for orthogonal model matrices ($\mathbf{W}= \mathbf{I}$)}
%\FigRule
%\vspace{-5ex}
\label{Orthhist}
\end{figure}
%------------------------------
In Figure-\ref{Orthhist}, we show the normalized histogram (normalized to make the total area one) of the components of the estimate vector
$\xhats_k, k = 1,2,3,4$ for the case of orthogonal models and compare it with the theoretical expression for the pdf obtained in (\ref{One/p1/Eqn/Ea13}).
%------------------------------
\begin{figure}[!tbp]
  \centering
  \begin{minipage}[b]{0.4\textwidth}
	\centering
    \includegraphics[width=4cm,height=5cm]{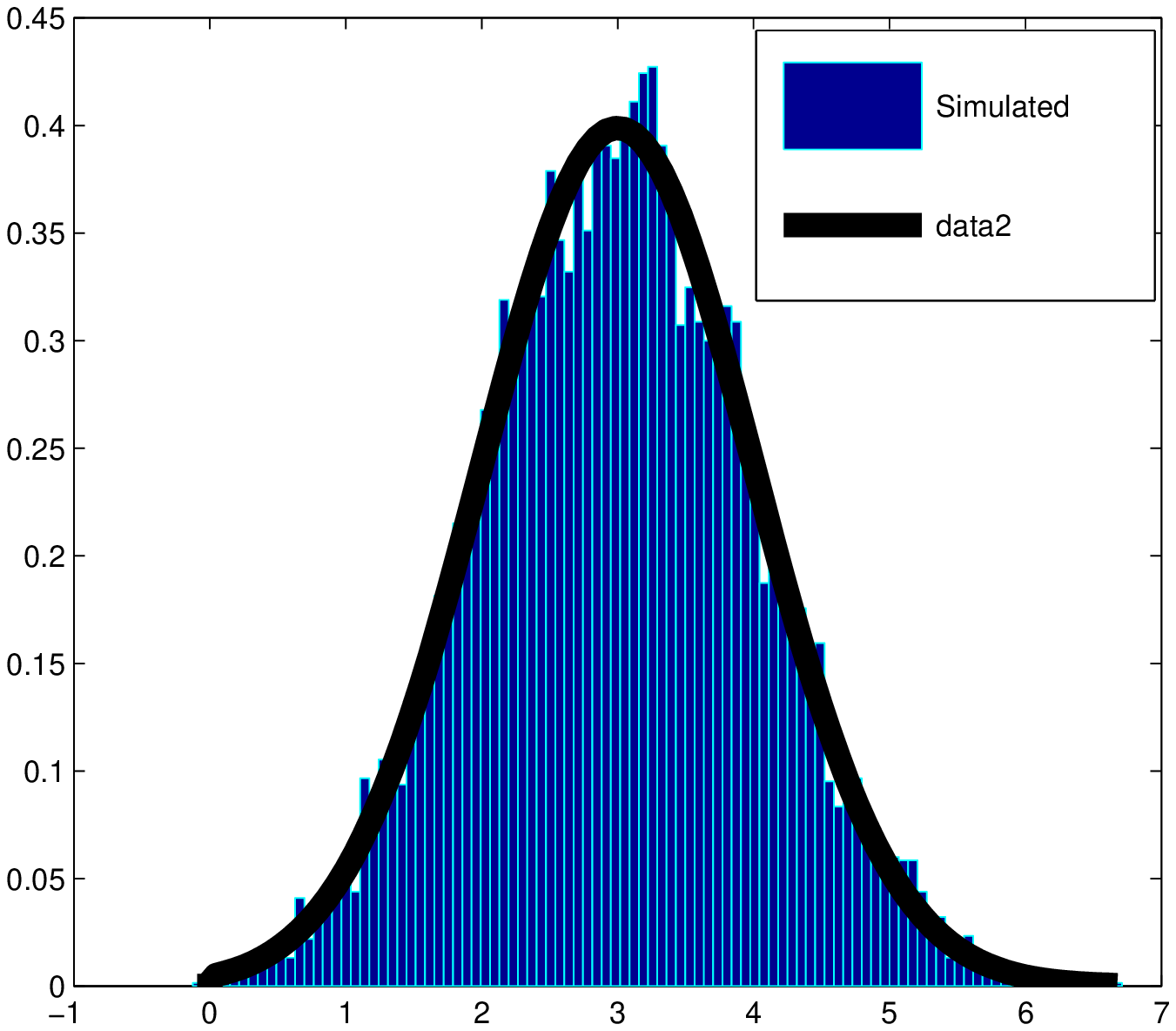}
    \caption{pdf of the $\zhats_{3}$ for full rank $W$}
    \label{Fullrank}
  \end{minipage}
  \hfill
  \begin{minipage}[b]{0.4\textwidth}
  	  \centering
    \includegraphics[width=4cm,height=5cm]{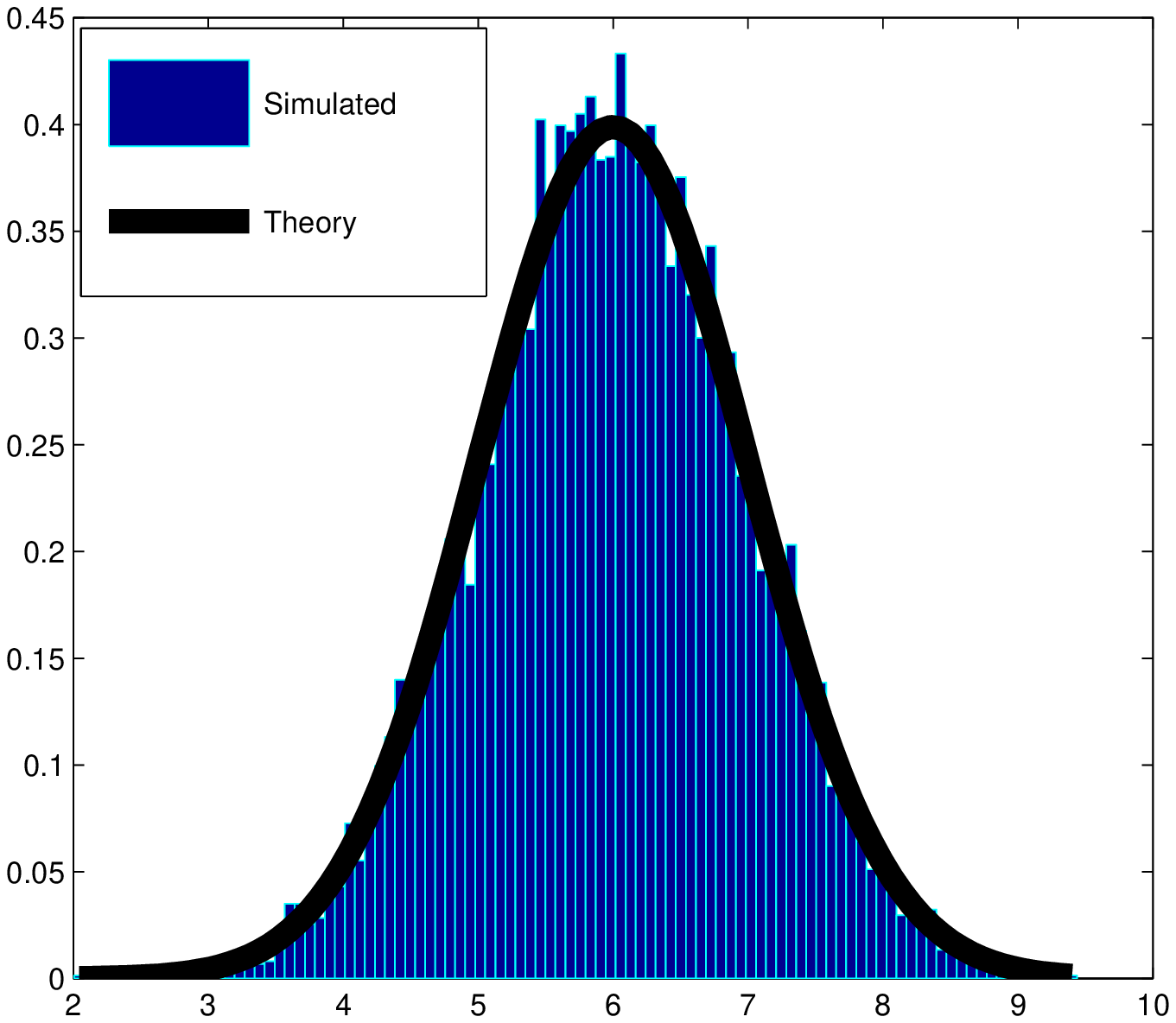}
    \caption{pdf of $\zhats_{5}$ for singular $\mathbf{W}$}
    \label{nonFullrank}
  \end{minipage}
\end{figure}
%-----------------------------
%------------------------------
In Figure-\ref{Fullrank}, we show the normalized histogram (normalized to make the total area one) of the components of the estimate vector
$\zhats_k, k = 3$ for the case of full rank model and compare it with the theoretical expression for the pdf obtained in (\ref{One/p1/Eqn/Ea16}). The simulated pdf does not match with (\ref{One/p1/Eqn/Ea16}) for the other components of $\zhat$.

In Figure-\ref{nonFullrank}, we show the normalized histogram (normalized to make the total area one) of the components of the estimate vector
$\zhats_k, k = 3$ for the case of full rank model and compare it with the theoretical expression for the pdf obtained in (\ref{One/p1/Eqn/Ea16c}). The simulated pdf does not match with (\ref{One/p1/Eqn/Ea16c}) for the other components of $\zhat$.

We observe from the figures that the theoretical pdfs follow the simulated  pdfs closely for all the components in case of orthogonal models and for $k$ corresponding to non-zero $\x$ in case of the other models.

Although, we have shown the simulations only for single source scenarios ($K=1$), the same simulations work for multiple source ($K>1$) scenarios and also in case of strong source-weak source scenarios, i.e. when $x_{max}>>x_{min}$.
%-------------------------------------------------------
\section{Concluding Remarks}\label{concl}
%-------------------------------------------------------
In this work, we have derived the cf of the LASSO estimator with an aim to get some insight on the distributional properties of the LASSO estimator. The expression of the cf contains derived in (\ref{One/p1/Eqn/Ea10}) contains $2^N$ summations, and hence does not simplify the evaluation of the distribution of the LASSO estimator except for the special case of orthogonal model matrix ($\W = \I$). So we use the Fourier-Slice theorem to calculate the one-dimensional projections of a linear transformation of the LASSO estimator. The approximate pdf of the one dimensional components this linear transformation is then found from the  projections by using the inversion theorem.

It is well-known from the Cramer Wold theorem that the one dimensional projections of a distribution are sufficient for the evaluation of the overall cf and hence the joint pdf of the LASSO estimator can be evaluated by its one dimensional projections, which is an interesting future work. As an application in statistical detection theory, the distribution of the estimator or (any of its functions, or test statistics) plays an important role to make decisions based on hypothesis. Hence, it may be an interesting future work to use the pdfs of the one-dimensional projections to propose test statistics for hypothesis testing. Many engineering applications like wireless communications and signal processing work with complex measurement models. Hence, a generalization of the distribution function for complex measurement models can be another interesting future work.
%---------------------------------------------
\appendix
%-------------------------------------------------------
\section{Proofs}\label{proofs}
In this section, we present the proofs of the theorems and corollaries discussed in Section-\ref{results}.
%-------------------------------------------------------
\subsection{Proof of Theorem-\ref{One/p1/Thm/Thm4}}\label{prf1}
%-------------------------------------------------------
The solution to the problem in (\ref{One/p1/Eqn/Ea4}) can be obtained in a straightforward manner by using the KKT conditions which results in the following implicit relation for $\xhat$.
%-----------------------------------
%------------------------------------
\begin{EqA}{One/p1/Eqn/Ea18}{rcl}
%-----------------------------------
\W\xhat + \tau\boldsymbol{\gamma} &=& \A^{T}\b = \btil\\
\gamma_k &\in &\left\{
%-----------------------------------
\begin{array}{rl}
\sgn(\xhats_k) & \text{if } \xhats_k \neq 0,\\
\left[-1,1\right] & \text{if } \xhats_k = 0
%-----------------------------------
\end{array} \right.
%-----------------------------------
\end{EqA}
%-----------------------------------
%-----------------------------------
where $\W =[\w_1,\w_2,\ldots,\w_{N}]= \A^{T}\A$. Let $\cJ = \{1,2,\ldots,N\}$ and $\eJ = \{k\in\cJ: |\gamma_k| = 1\}$, then any LASSO solution $\xhat$ satisfies
%-----------------------------------
%------------------------------------
\begin{EqA}{One/p1/Eqn/Ea19}{rcll}
%-----------------------------------
\xhat_{\cJ\setminus\eJ} = \mathbf{0};&\W_{\eJ}\xhat_{\eJ} + \tau\gm_{\eJ} &=& \A^{T}_{\eJ}\b = \br
%-----------------------------------
\end{EqA}
%-----------------------------------
%-----------------------------------
Our aim is now to calculate the cf of $\xhat$ using (\ref{One/p1/Eqn/Ea19}). Since, $\xhat_{\cJ\setminus\eJ} = \mathbf{0}$, it does not contribute in the evaluation of the cf. Hence, it is enough to consider the first part of (\ref{One/p1/Eqn/Ea18}).

Now, we first evaluate the cf of right hand side (R.H.S) of (\ref{One/p1/Eqn/Ea18}). We observe that $\btil$ is a multivariate Gaussian random variable with mean $\W\x$ and variance $\sigma^2\W$. Hence, the cf of $\btil$, by definition is
%-----------------------------------
%------------------------------------
\begin{EqA}{One/p1/Eqn/Ea20}{rcl}
%-----------------------------------
\cf_{\btil}(\u) &=& \exp\big(\iu\u^{T}\W\x-\frac{\sigma^2}{2}\u^{T}\W\u\big)
%-----------------------------------
\end{EqA}
%-----------------------------------
%-----------------------------------
Now, we need to calculate the cf of the left hand side (L.H.S) of (\ref{One/p1/Eqn/Ea18}). We first define $c_j = \u^{T}\w_j$, $\xhat^{-} = [\xhats_2,\xhats_3,\ldots,\xhats_N]$ and d$\xhat^{-}=$d$\xhats_2$d$\xhats_2\ldots$d$\xhats_N$. We have,
%-----------------------------------
%------------------------------------
\begin{EqA}{One/p1/Eqn/Ea21}{rcl}
%-----------------------------------
&&\cf_{\mathsmaller{\W\xhat + \tau\sgn(\xhat)}}(\u) = \mean_{\mathsmaller{\W\xhat + \tau\sgn(\xhat)}}\bigg\{\exp\bigg(\iu\u^{T}\big(\W\xhat + \tau\sgn(\xhat)\big)\bigg)\bigg\}\\
%-----------------------------------
&=&\mean_{\xhat}\bigg\{\exp\bigg(\iu\SUM{j}{1}{N}{\bigg(\xhats_{j}\u^{T}\w_{j}+\tau u_{j}\sgn(\xhats_{j})\bigg)}\bigg)\bigg\}\\
%-----------------------------------
&=&\INT{\re^{N}}{}{\hspace{1mm}f_{\mathsmaller{\xhat}}(\xhat)\exp\bigg(\iu\SUM{j}{1}{N}{\bigg(\xhats_{j}\u^{T}\w_{j}+\tau u_{j}\sgn(\xhats_{j})\bigg)}\bigg)\mathrm{d}\xhat}\\
%-----------------------------------
&=&\INT{\re^{N}}{}{\hspace{1mm}f_{\mathsmaller{\xhat}}(\xhats_1,\ldots,\xhats_N)\exp\bigg(\iu\SUM{j}{1}{N}{\bigg(\xhats_{j}\u^{T}\w_{j}+\tau u_{j}\sgn(\xhats_{j})\bigg)}\bigg)\mathrm{d}\xhat}\\
%-----------------------------------
&=&\INT{\re^{N}}{}{\hspace{1mm}f_{\mathsmaller{\xhat}}(\xhats_1,\ldots,\xhats_N)\exp\bigg(\iu\SUM{j}{1}{N}{\bigg(\xhats_{j}c_{j}+\tau u_{j}\sgn(\xhats_{j})\bigg)}\bigg)\mathrm{d}\xhat}\\
%-----------------------------------
&=&\INT{\re^{N-1}}{}{\Bigg(\INT{\re}{}{f_{\mathsmaller{\xhat}}(\xhat)\exp(\iu(\xhats_{1}c_{1}+ \tau u_{1}\sgn(\xhats_{1})))\mathrm{d}\xhats_{1}\Bigg)}\exp\bigg(\iu\SUM{j}{2}{N}{\bigg(\xhats_{j}c_{j}+ \tau u_{j}\sgn(\xhats_{j})\bigg)}\bigg) \mathrm{d}\xhat^{-}}\\
%-----------------------------------
&=&\INT{\re^{N-1}}{}{\mathcal{J}(c_1,u_1,\xhat^{-})\exp\bigg(\iu\SUM{j}{2}{N}{\bigg(\xhats_{j}c_{j}+ \tau u_{j}\sgn(\xhats_{j})\bigg)}\bigg) \mathrm{d}\xhat^{-}}
%-----------------------------------
\end{EqA}
%-----------------------------------
%-----------------------------------
We first evaluate $\mathcal{J}(c_1,u_1,\xhat^{-})$. Defining $H(x)$ as the Heaviside step function, we have
%-----------------------------------
%-----------------------------------
\begin{eqnarray}
%-----------------------------------
&&\mathcal{J}(c_1,u_1,\xhat^{-})=\INT{-\infty}{\infty}{f_{\mathsmaller{\xhat}}(\xhat)\exp(\iu(\xhats_{1}c_{1}+ \tau u_{1}\sgn(\xhats_{1})))\mathrm{d}\xhats_{1}} \nonumber\\
%-----------------------------------
&&=\INT{-\infty}{0}{f_{\mathsmaller{\xhat}}(\xhat)e^{\iu(\xhats_{1}c_{1}- \tau u_{1})}\mathrm{d}\xhats_{1}}+ \INT{0}{\infty}{f_{\mathsmaller{\xhat}}(\xhat)e^{\iu(\xhats_{1}c_{1}+\tau u_{1})}\mathrm{d}\xhats_{1}}\nonumber\\
%-----------------------------------
&&=e^{-\iu\tau u_{1}}\INT{-\infty}{0}{f_{\mathsmaller{\xhat}}(\xhat)e^{\iu(\xhats_{1}c_{1})}\mathrm{d}\xhats_{1}}+
e^{\iu\tau u_{1}}\INT{0}{\infty}{f_{\mathsmaller{\xhat}}(\xhat)e^{\iu(\xhats_{1}c_{1})}\mathrm{d}\xhats_{1}}\nonumber\\
%-----------------------------------
&&=e^{-\iu\tau u_{1}}\INT{-\infty}{\infty}{f_{\mathsmaller{\xhat}}(\xhat)H(-\xhats_{1})e^{\iu(\xhats_{1}c_{1})}\mathrm{d}\xhats_{1}}+e^{\iu\tau u_{1}}\INT{-\infty}{\infty}{f_{\mathsmaller{\xhat}}(\xhat)H(\xhats_{1})e^{\iu(\xhats_{1}c_{1})}\mathrm{d}\xhats_{1}}
\label{One/p1/Eqn/Ea22-4}\\
%-----------------------------------
&&=e^{-\iu\tau u_{1}}\bigg(\ocf_{\xhats_1}(c_{1},\xhat^{-})\star \cf_{H}(-c_{1})\bigg)+e^{\iu\tau u_{1}}\bigg(\ocf_{\xhats_1}(c_{1},\xhat^{-})\star \cf_{H}(c_{1})\bigg)\label{One/p1/Eqn/Ea22-5}\\
%-----------------------------------
&&=e^{-\iu\tau u_{1}}\bigg(\ocf_{\xhats_{1}}(c_{1},\xhat^{-})\star\bigg(\frac{1}{2\pi jc_{1}}+\frac{\delta(c_{1})} {2}\bigg)\bigg)+e^{\iu\tau u_{1}}\bigg(\ocf_{\xhats_{1}}(c_{1},\xhat^{-})\star\bigg(\frac{-1}{2\pi jc_{1}}+ \frac{\delta(c_{1})}{2}\bigg)\bigg)\label{One/p1/Eqn/Ea22-6}\\
%-----------------------------------
&&=\ocf_{\xhats_{1}}(c_{1},\xhat^{-})\cos(\tau u_1)+\bigg(\ocf_{\xhats_{1}}(c_{1},\xhat^{-})\star\Big(\frac{-1}{\pi c_{1}}\Big)\bigg)\sin(\tau u_1)\nonumber\\
%-----------------------------------
&&=\ocf_{\xhats_{1}}(c_{1},\xhat^{-})\cos(\tau u_1)+\ocf_{\xhats_{1}}(\overline{c}_{1},\xhat^{-})\sin(\tau u_1)\label{One/p1/Eqn/Ea22-8}
%-----------------------------------
\end{eqnarray}
%----------------------------------
%-----------------------------------
Here, $\ocf_{\xhats_{1}}(c_{1},\xhat^{-})$ denotes the one-dimensional Fourier transform along $c_1$ and $\ocf_{\xhats_{1}}(\overline{c}_{1},\xhat^{-})= \bigg(\ocf_{\xhats_{1}}(c_{1},\xhat^{-})\star\Big(\frac{-1}{\pi c_{1}}\Big)\bigg)$, denotes the one-dimensional Hilbert transform along $c_1$. Similarly, the Fourier and Hilbert transforms over each dimension is absorbed into the notation of the cf. We have used the Fourier-transform convolution theorem to obtain (\ref{One/p1/Eqn/Ea22-5}) from (\ref{One/p1/Eqn/Ea22-4}) and we have used the Fourier transform of the Heaviside step function in (\ref{One/p1/Eqn/Ea22-6}).

Now using the value of $\mathcal{J}(c_1,u_1,\xhat^{-})$ obtained from (\ref{One/p1/Eqn/Ea22-8}) in (\ref{One/p1/Eqn/Ea21}), we have
%-----------------------------------
%------------------------------------
\begin{EqA}{One/p1/Eqn/Ea23}{rcl}
%-----------------------------------
\cf_{\mathsmaller{\W\xhat + \tau\sgn(\xhat)}}(\u)&=&\INT{\re^{N-1}}{}{\mathcal{J}(c_1,u_1,\xhat^{-}) \exp\bigg(\iu\SUM{j}{2}{N}{\bigg(\xhats_{j}c_{j}+ \tau u_{j}\sgn(\xhats_{j})\bigg)}\bigg) \mathrm{d}\xhat^{-}}\\
&=&\cos(\tau u_1)\INT{\re^{N-1}}{}{\cf_{\xhats_{1}}(c_{1},\xhat^{-})\exp\bigg(\iu\SUM{j}{2}{N}{\bigg(\xhats_{j}c_{j}+\tau u_{j}\sgn(\xhats_{j})\bigg)}\bigg) \mathrm{d}\xhat^{-}}\\
&&+\sin(\tau u_1)\INT{\re^{N-1}}{} {\cf_{\xhats_{1}}(\overline{c}_{1},\xhat^{-})\exp\bigg(\iu\SUM{j}{2}{N}{\bigg(\xhats_{j}c_{j}+\tau u_{j}\sgn(\xhats_{j})\bigg)}\bigg) \mathrm{d}\xhat^{-}}
%-----------------------------------
\end{EqA}
%-----------------------------------
%-----------------------------------
Now, integrating over $\xhats_2,\xhats_3,\ldots\xhats_N$ in the same manner, we obtain
%-----------------------------------
%------------------------------------
\begin{EqA}{One/p1/Eqn/Ea24}{rl}
%-----------------------------------
\cf_{\mathsmaller{\W\xhat + \tau\sgn(\xhat)}}(\u)=&\cos(\tau u_1)\cos(\tau u_2)\ldots\cos(\tau u_N)\cf_{\xhat}(c_1,\ldots,c_N)\\
&+\sin(\tau u_1)\cos(\tau u_2)\ldots\cos(\tau u_N)\cf_{\xhat}(\overline{c}_1,c_2,\ldots,c_N)\\
&+\cos(\tau u_1)\sin(\tau u_2)\ldots\cos(\tau u_N)\cf_{\xhat}(c_1,\overline{c}_2,\ldots,c_N)\\
&\ldots+(\text{All combinations of $1$ sine term and $N-1$ cosine terms})\\
&\ldots+\sin(\tau u_1)\sin(\tau u_2)\cos(\tau u_3)\ldots\cos(\tau u_N)\cf_{\xhat}(\overline{c}_1,\overline{c}_2,c_3,\ldots,c_N)\\
&\ldots+\cos(\tau u_1)\sin(\tau u_2)\sin(\tau u_3)\ldots\cos(\tau u_N)\cf_{\xhat}(c_1,\overline{c}_2,\overline{c}_3,c_4,\ldots,c_N)\\
&\ldots+\text{(All combinations of $2$ sine terms and $N-2$ cosine terms)}\\
&\hspace{5 cm}\vdots\\
&\ldots+\text{(All combinations of $k$ sine terms and $N-k$ cosine terms)}\\
&\hspace{5 cm}\vdots\\
&\ldots+\sin(\tau u_1)\ldots+\sin(\tau u_N)\cf_{\xhat}(\overline{c}_1,\ldots,\overline{c}_N)
%-----------------------------------
\end{EqA}
%-----------------------------------
%-----------------------------------
which is equal to expression (\ref{One/p1/Eqn/Ea10})
%-------------------------------------------------------
\subsection{Proof of Corollary-\ref{One/p1/Thm/Thm5}}\label{prf2}
%-------------------------------------------------------
%-----------------------------------
If $\W = \I$, where $\I$ is the identity matrix, or any diagonal matrix $\D$ (in general), then the estimates $\hat{x}_k, k=1,2,\ldots,N$ are all independent and hence the cf of the LASSO estimator is given by,
%-----------------------------------
%------------------------------------
\begin{EqA}{One/p1/Eqn/Ea25}{rcl}
%-----------------------------------
\PROD{k}{1}{N}{\hspace{1mm}\Bigg(\cf_{\xhats_{k}}(u_{k})\cos(\tau u_{k})+ \cf_{\xhats_{k}}(\overline{u}_{k})\sin(\tau u_{k})\Bigg)} &=&\PROD{k}{1}{N}{\hspace{1mm}\exp\big(\iu u_{k}x_{k}-\frac{\sigma^2}{2}u^{2}_{k}\big)}
%-----------------------------------
\end{EqA}
%-----------------------------------
%-----------------------------------
which can be expressed as,
%-----------------------------------
%------------------------------------
\begin{EqA}{One/p1/Eqn/Ea12app}{rcl}
%-----------------------------------
\cf_{\xhats_{k}}(u_{k})\cos(\tau u_{k})+ \cf_{\xhats_{k}}(\overline{u}_{k})\sin(\tau u_{k})&=& \exp\big(\iu u_{k}x_{k}-\frac{\sigma^2}{2}u^{2}_{k}\big), k = 1,2,\ldots,N
%-----------------------------------
\end{EqA}
%-----------------------------------
%-----------------------------------
Hence, the marginal probability distribution functions of the individual components, $\hat{x}_k, k=1,2,\ldots,N$ is given by simply applying the inversion theorem to the (\ref{One/p1/Eqn/Ea12app}) as,
%-----------------------------------
%------------------------------------
\begin{EqA}{One/p1/Eqn/Ea26}{rcl}
%-----------------------------------
\frac{1}{\sqrt{2\pi\sigma^2}}\exp\big(-\frac{(\xhats_i-x_k)^2}{2\sigma^2}\big)&=&\frac{1}{2}f_{\xhats_k}(\xhats_k)\star\big[\delta(\xhats_k-\tau)+
\delta(\xhats_k+\tau)\big]\\
&&+\frac{1}{2}\sgn(\xhats_k)f_{\xhats_k}(\xhats_k)\star \big[\delta(\xhats_k-\tau)-\delta(\xhats_k+\tau)\big]\\
&=&\frac{1}{2}\big[f_{\xhats_k}(\xhats_k-\tau)+f_{\xhats_k}(\xhats_k+\tau)\big] +\frac{1}{2}\sgn(\xhats_k-\tau) f_{\xhats_k}(\xhats_k-\tau)\\
&&-\frac{1}{2}\sgn(\xhats_k+\tau)f_{\xhats_k}(\xhats_k+\tau)
%-----------------------------------
\end{EqA}
%-----------------------------------
which gives,
%------------------------------------
\begin{EqA}{One/SetA1/p1/Eqn/Ea26-a}{rcl}
%-----------------------------------
\frac{1}{\sqrt{2\pi\sigma^2}}\exp\big(-\frac{(\xhats_k-x_k)^2}{2\sigma^2}\big)& = & \left\{
%-----------------------------------
\begin{array}{rl}
f_{\xhats_k}(\xhats_k-\tau)& \text{if } \xhats_k >\tau,\\
f_{\xhats_k}(\xhats_k+\tau) & \text{if } \xhats_k <-\tau
%-----------------------------------
\end{array} \right.
%-----------------------------------
\end{EqA}
%-----------------------------------
%-----------------------------------
Hence,
%-----------------------------------
%------------------------------------
\begin{EqA}{One/p1/Eqn/Ea13app}{rcl}
%-----------------------------------
f_{\xhats_k}(\xhats_k)& = & \left\{
%-----------------------------------
\begin{array}{rl}
\frac{1}{\sqrt{2\pi\sigma^2}}\exp\big(-\frac{(\xhats_k+\tau-x_k)^2}{2\sigma^2}\big)& \text{if } \xhats_k>0,\\
\frac{1}{\sqrt{2\pi\sigma^2}}\exp\big(-\frac{(\xhats_k-\tau-x_k)^2}{2\sigma^2}\big) & \text{if } \xhats_k<0.
%-----------------------------------
\end{array} \right.
%-----------------------------------
\end{EqA}
%-------------------------------------------------------
\subsection{Proof of Corollary-\ref{One/p1/Thm/Thm6}}\label{prf3}
%-------------------------------------------------------
%-----------------------------------
If $\W$ is any full rank matrix, we first perform slicing by substituting $\u = u\e_{k}$ in equation (\ref{One/p1/Eqn/Ea24}) to obtain
%-----------------------------------
%------------------------------------
\begin{EqA}{One/p1/Eqn/Ea27}{rcl}
%-----------------------------------
\cf_{\mathsmaller{\W\xhat + \tau\sgn(\xhat)}}(u)&=&\cos(\tau u)\mathcal{S}_{\u}[\cf_{\xhat}(\W\u)] + \sin(\tau u) \mathcal{S}_{\u}[\cf_{\xhat}(\W\u)\star(\frac{-1}{\pi c_k})]
%-----------------------------------
\end{EqA}
%-----------------------------------
%-----------------------------------
Now, we define $\zhat = \W\xhat$ and $\xhat = \H\zhat$, where $\H = \W^{-1}$. Let $\h_k$ be the $k^{th}$ column of $\H$, we have $f(\zhat) = |\W|^{-1} f(\xhat)$ and $\mathcal{S}_{\u}[\cf_{\xhat}(\W\u)]= \mathcal{S}_{\u}\odot\W^{-1}\odot\cf_{\xhat} = \mathcal{S}_{\u}\odot\cf_{\zhat}\odot\W = \cf_{\zhats_k}\odot\mathcal{I}_{\zhats_k}\odot\W = \cf_{\zhats_k}(u)$, where $\zhats_k$ is the $k^{th}$ component of $\zhat$ and we obtain the last equation by applying the generalized Fourier slice theorem (Theorem-\ref{One/p1/Thm/Thm3}). We evaluate the term, $\mathcal{S}_{\u}[\cf_{\xhat}(\W\u)\star(\frac{-1}{\pi c_k})]$ below by noting that
$\frac{-1}{\pi u_k}$ is the Fourier transform of $\iu\sgn(\xhats_k)$,
%-----------------------------------
%------------------------------------
\begin{EqA}{One/p1/Eqn/Ea28-1}{rcl}
%-----------------------------------
\cf_{\xhat}(\W\u)\star(\frac{-1}{\pi c_k})&=& \iu\INT{\re^{N}}{}{\hspace{1mm}f_{\mathsmaller{\xhat}}(\xhat)\sgn(x_k)\exp(\iu\u^{T}\W\xhat)\mathrm{d}\xhat}\\
%-----------------------------------
&=& \iu\INT{\re^{N}}{}{\hspace{1mm}f_{\mathsmaller{\zhat}}(\zhat)\sgn(\h^{T}\zhat)\exp(\iu\u^{T}\zhat)\mathrm{d}\zhat}\\
%-----------------------------------
\mathcal{S}_{\u}[\cf_{\xhat}(\W\u)\star(\frac{-1}{\pi c_k})]&=& \iu\INT{\re^{N}}{}{\hspace{1mm}f_{\mathsmaller{\zhat}}(\zhat)\sgn(\h_{k}^{T}\zhat)\exp(\iu u\zhats_k)\mathrm{d}\zhat}
%-----------------------------------
\end{EqA}
%-----------------------------------
%-----------------------------------
We now make an approximation that $\sgn(\h_{k}^{T}\zhat)\approx \sgn(\zhats_k)$ for $k$ corresponding to large $|\xhats_{k}|$ as explained in Theorem-\ref{One/p1/Thm/Thm6}. Hence, the term $\mathcal{S}_{\u}[\cf_{\xhat}(\W\u)\star(\frac{-1}{\pi c_k})] \approx\cf_{\zhats_k}(\overline{u})$. Substituting for $\mathcal{S}_{\u}[\cf_{\xhat}(\W\u)]$ and the approximated expression of $\mathcal{S}_{\u}[\cf_{\xhat}(\W\u)\star(\frac{-1}{\pi c_k})]$ for $k$ corresponding to large $|\xhats_{k}|$ and equating the above expression to the slice (w.r.t to $u$) of the cf of $\btil$, we have
%-----------------------------------
%------------------------------------
\begin{EqA}{One/p1/Eqn/Ea16app}{rcl}
%-----------------------------------
\cos(\tau u)\cf_{\zhats_k}(u)+ \sin(\tau u)\cf_{\zhats_k}(\overline{u})
&=& \exp(-u^{2}\frac{\sigma^2}{2}w_{kk})\exp\Big(ju(\w_{k}^{T}\x)\Big)
%-----------------------------------
\end{EqA}
%-----------------------------------
%-----------------------------------
Hence, the marginal pdf of the individual components, $\hat{z}_k$ for $k$ corresponding to large $|\xhats_{k}|$ is given by simply applying the inversion theorem to the (\ref{One/p1/Eqn/Ea16app}) as,
%-----------------------------------
%------------------------------------
\begin{EqA}{One/p1/Eqn/Ea17}{rcl}
%-----------------------------------
f_{\zhats_k}(\zhats_k)& = & \left\{
%-----------------------------------
\begin{array}{rl}
\frac{1}{\sqrt{2\pi\sigma^2w_{kk}}}\exp\big(-\frac{(\zhats_k+\tau-\w^{T}_{k}\x)^2}{2\sigma^2}\big)& \text{if } \zhats_k>0,\\
\frac{1}{\sqrt{2\pi\sigma^2w_{kk}}}\exp\big(-\frac{(\zhats_k-\tau-\w^{T}_{k}\x)^2}{2\sigma^2}\big) & \text{if } \zhats_k<0.
%-----------------------------------
\end{array} \right.
%-----------------------------------
\end{EqA}
%-----------------------------------
%-------------------------------------------------------
\subsection{Proof of Corollary-\ref{One/p1/Thm/Thm7}}\label{prf4}
%-------------------------------------------------------
%-----------------------------------
For any general $\W$, we again perform slicing by substituting $\u = u\e_{k}$ in equation (\ref{One/p1/Eqn/Ea24}) to obtain (\ref{One/p1/Eqn/Ea27}). Again, we define $\zhat = \W\xhat$, $\xhat = \H\zhat$, where $\H$ is now $\H = \W^{\dagger}$ and let $\h_k$ be the $k^{th}$ column of $\H$.

Now, as in Proof-\ref{prf3}, we have $\mathcal{S}_{\u}[\cf_{\xhat}(\W\u)]= \mathcal{S}_{\u}\odot\cf_{\zhat}\odot\W = \cf_{\zhats_k}\odot\mathcal{I}_{\zhats_k}\odot\W = \cf_{\zhats_k}(u)$, where $\zhats_k$ is the $k^{th}$ element of $\zhat$ and we obtain the last equation by applying the generalized Fourier slice theorem (Theorem-\ref{One/p1/Thm/Thm3}). The term, $\mathcal{S}_{\u}[\cf_{\xhat}(\W\u)\star(\frac{-1}{\pi c_k})]$ is again equal to,
%-----------------------------------
%------------------------------------
\begin{EqA}{One/p1/Eqn/Ea15}{rcl}
%-----------------------------------
\mathcal{S}_{\u}[\cf_{\xhat}(\W\u)\star(\frac{-1}{\pi c_k})]&=&\iu \INT{\re^{N}}{}{\hspace{1mm}f_{\mathsmaller{\zhat}}(\zhat)\sgn(\h_{k}^{T}\zhat)\exp(ju\zhats_k)\mathrm{d}\zhat}
%-----------------------------------
\end{EqA}
%-----------------------------------
%-----------------------------------
We again make the approximation $\sgn(\h_{k}^{T}\zhat)\approx \sgn(\zhats_k)$ for $k$ corresponding to large $|\xhats_{k}|$. Substituting for
$\mathcal{S}_{\u}[\cf_{\xhat}(\W\u)]$ and the approximated expression of $\mathcal{S}_{\u}[\cf_{\xhat}(\W\u)\star(\frac{-1}{\pi c_k})]$ for $k$ corresponding to large $|\xhats_{k}|$ in equation (\ref{One/p1/Eqn/Ea20}) and equating the above expression to the slice (w.r.t to $u$) of the cf of $\btil$, we have
%-----------------------------------
%------------------------------------
\begin{EqA}{One/p1/Eqn/Ea16a}{rcl}
%-----------------------------------
\cos(\tau u)\cf_{\zhats_k}(u)+ \sin(\tau u)\cf_{\zhats_k}(\overline{u})
&=& \exp(-u^{2}\frac{\sigma^2}{2}w_{kk})\exp\Big(ju(\w_{k}^{T}\x)\Big)
%-----------------------------------
\end{EqA}
%-----------------------------------
%-----------------------------------
Hence, the marginal pdf of the individual components, $\hat{z}_k$ for $k$ corresponding to large $|\xhats_{k}|$ is given by simply applying the inversion theorem to the (\ref{One/p1/Eqn/Ea16a}) as,
%-----------------------------------
%------------------------------------
\begin{EqA}{One/p1/Eqn/Ea17}{rcl}
%-----------------------------------
f_{\zhats_k}(\zhats_k)& = & \left\{
%-----------------------------------
\begin{array}{rl}
\frac{1}{\sqrt{2\pi\sigma^2w_{kk}}}\exp\big(-\frac{(\zhats_k+\tau-\w^{T}_{k}\x)^2}{2\sigma^2}\big)& \text{if } \zhats_k>0,\\
\frac{1}{\sqrt{2\pi\sigma^2w_{kk}}}\exp\big(-\frac{(\zhats_k-\tau-\w^{T}_{k}\x)^2}{2\sigma^2}\big) & \text{if } \zhats_k<0.
%-----------------------------------
\end{array} \right.
%-----------------------------------
\end{EqA}
%-----------------------------------
%-------------------------------------------------------
\subsection{Evaluation of $\mathcal{S}_{\u}[\cf_{\xhat}(\W\u)\star(\frac{-1}{\pi c_k})]$}\label{hilbertval}
%-------------------------------------------------------
In this section, we evaluate $\mathcal{S}_{\u}[\cf_{\xhat}(\W\u)\star(\frac{-1}{\pi c_k})]$ with the assumption that $\zhat$ has a multivariate Gaussian distribution of with mean $\m$ and co-variance matrix $\R$. We use $\phi(\zhat_{\mathsmaller{N}},\m_{\mathsmaller{N}}, \R_{\mathsmaller{N}})$ to denote that the random vector $\zhat_{\mathsmaller{N}}$ of length $N$ has a multivariate Gaussian distribution with mean vector $\m_{\mathsmaller{N}}$ of length $N$ and co-variance matrix $\R_{\mathsmaller{N}}$ of size $N\times N$ and $\Phi$ is used to denote the cumulative distribution function (cdf) of the normal distribution. Let $\H = \W^{\dagger}$ and $\h_{k}$ as the $k^{th}$ column of $\H$. We have,
%-----------------------------------
%------------------------------------
\begin{EqA}{One/p1/Eqn/Ea29}{rcl}
%-----------------------------------
\mathcal{S}_{\u}[\cf_{\xhat}(\W\u)\star(\frac{-1}{\pi c_k})]&=&\iu\INT{\re^{N}}{}{\hspace{1mm}f_{\mathsmaller{\zhat}}(\zhat)\exp(ju\zhats_k)\sgn(\h_{k}^{T}\zhat)\mathrm{d}\zhat}
%-----------------------------------
\end{EqA}
%-----------------------------------
%-----------------------------------
Without loss of generality, we choose $k = 1$. We partition
$\m_{\mathsmaller{N}}$, $\zhat_{\mathsmaller{N}}$ and
$\R_{\mathsmaller{N}}$ as $\m = [\m_{\mathsmaller{N-1}}, m_{\mathsmaller{N}}]^{T}$, $\zhat = [\zhat_{\mathsmaller{N-1}},\zhats_{\mathsmaller{N}}]^{T}$, $\R =\begin{bmatrix} \R_{\mathsmaller{N-1}}&\r_{\mathsmaller{N}}\\\r_{\mathsmaller{N}}^{T}&r_{\mathsmaller{NN}}\end{bmatrix}$ respectively. Let $\s = \frac{-1}{h_{\mathsmaller{1N}}}[h_{\mathsmaller{11}},h_{\mathsmaller{12}},\ldots, h_{\mathsmaller{1N}}]^{T}$, then $\sgn(\h_{\mathsmaller{1}}^{T}\zhat) =\mp1$ depending on $\zhats_{\mathsmaller{N}}\lessgtr\s^{T}\zhat_{\mathsmaller{N-1}}$. We have,
%-----------------------------------
%------------------------------------
\begin{EqA}{One/p1/Eqn/Ea30-1}{rcl}
%-----------------------------------
\mathcal{S}_{\u}[\cf_{\xhat}(\W\u)\star(\frac{-1}{\pi c_k})]&=&\iu\INT{\re^{N}}{}{\hspace{1mm}f_{\mathsmaller{\zhat}}(\zhat)\exp(ju\zhats_1)\sgn(\h_{1}^{T}\zhat)\mathrm{d}\zhat}\\
%-----------------------------------
&=&\iu\INT{\re^{N-1}}{}{\Bigg(-\INT{-\infty}{\s^{T}\zhat_{\mathsmaller{N-1}}}{f_{\mathsmaller{\zhat}}(\zhat)\mathrm{d}\zhats_{\mathsmaller{N}}}+\INT{\s^{T}\zhat_{\mathsmaller{N-1}}}{\infty}{f_{\mathsmaller{\zhat}}(\zhat)\mathrm{d}\zhats_{\mathsmaller{N}}}\Bigg)\exp(ju\zhats_1)\mathrm{d}\zhat_{\mathsmaller{N-1}}}\\
%-----------------------------------
&=&\iu\INT{\re^{N-1}}{}{\Bigg(-2\INT{-\infty}{\s^{T}\zhat_{\mathsmaller{N-1}}}{f_{\mathsmaller{\zhat}}(\zhat)\mathrm{d}\zhats_{\mathsmaller{N}}}+\INT{-\infty}{\infty}{f_{\mathsmaller{\zhat}}(\zhat)\mathrm{d}\zhats_{\mathsmaller{N}}}\Bigg)\exp(ju\zhats_1)\mathrm{d}\zhat_{\mathsmaller{N-1}}}\\
%-----------------------------------
&=&-2\iu\INT{\re^{N-1}}{}{\Bigg(\INT{-\infty}{\s^{T}\zhat_{\mathsmaller{N-1}}}{f_{\mathsmaller{\zhat}}(\zhat_{\mathsmaller{N-1}},\zhats_{\mathsmaller{N}})\mathrm{d}\zhats_{\mathsmaller{N}}}\Bigg)\exp(ju\zhats_1)\mathrm{d}\zhat_{\mathsmaller{N-1}}}+\iu\cf_{\zhats_{1}}(u)
%-----------------------------------
\end{EqA}
%-----------------------------------
%-----------------------------------
We now evaluate the first integral below.
%-----------------------------------
%------------------------------------
\begin{EqA}{One/p1/Eqn/Ea31}{rcl}
%-----------------------------------
J &=&-2\INT{\re^{N-1}}{}{\Bigg(\INT{-\infty}{\s^{T}\zhat_{\mathsmaller{N-1}}}{f_{\mathsmaller{\zhat}}(\zhat_{\mathsmaller{N-1}},\zhats_{\mathsmaller{N}})\mathrm{d}\zhats_{\mathsmaller{N}}}\Bigg)\exp(ju\zhats_1)\mathrm{d}\zhat_{\mathsmaller{N-1}}}
%-----------------------------------
\end{EqA}
%-----------------------------------
%-----------------------------------
Let us consider the inner integral, defining $p_{\mathsmaller{N}} = m_{\mathsmaller{N}}+\r^{T}_{\mathsmaller{N}}\R^{\dagger}_{\mathsmaller{N-1}}(\zhat_{\mathsmaller{N-1}}-\m_{\mathsmaller{N-1}})$ and $q_{\mathsmaller{N}} = r_{\mathsmaller{NN}}-\r^{T}_{\mathsmaller{N}}\R^{\dagger}_{\mathsmaller{N-1}}\r_{\mathsmaller{N}}$ and using the fact that $f(\zhat)$ is multivariate Gaussian, we have
%-----------------------------------
%------------------------------------
\begin{EqA}{One/p1/Eqn/Ea32}{rcl}
%-----------------------------------
\INT{-\infty}{\s^{T}\zhat_{\mathsmaller{N-1}}}{f_{\mathsmaller{\zhat}}(\zhat_{\mathsmaller{N-1}},\zhats_{\mathsmaller{N}})\mathrm{d}\zhats_{\mathsmaller{N}}}&=&\INT{-\infty}{\s^{T}\zhat_{\mathsmaller{N-1}}}{\phi(\zhat_{\mathsmaller{N-1}},\m_{\mathsmaller{N-1}},\R_{\mathsmaller{N-1}})N(\zhats_{\mathsmaller{N}},
p_{\mathsmaller{N}},q_{\mathsmaller{N}})\mathrm{d}\zhats_{\mathsmaller{N}}}
\\
%-----------------------------------
&=&\phi(\zhat_{\mathsmaller{N-1}},\m_{\mathsmaller{N-1}},\R_{\mathsmaller{N-1}})\Phi\Bigg(\frac{\s^{T}\zhat_{\mathsmaller{N-1}}-p_{\mathsmaller{N}}}{\sqrt{q_{\mathsmaller{N}}}}\Bigg)\\
%-----------------------------------
&=&\phi(\zhat_{\mathsmaller{N-1}},\m_{\mathsmaller{N-1}},\R_{\mathsmaller{N-1}})\Phi\Bigg(\frac{\g_{\mathsmaller{N-1}}^{T}\zhat_{\mathsmaller{N-1}}-k_{\mathsmaller{N-1}}}{
\sqrt{q_{\mathsmaller{N}}}}\Bigg)
%-----------------------------------
\end{EqA}
%-----------------------------------
%-----------------------------------
where $\g_{\mathsmaller{N-1}}^{T} = \s^{T}-\r^{T}_{\mathsmaller{N}}\R^{\dagger}_{\mathsmaller{N-1}}$ and $k_{\mathsmaller{N-1}} = \r^{T}_{\mathsmaller{N}}\R^{\dagger}_{\mathsmaller{N-1}}\m_{\mathsmaller{N-1}}-m_{\mathsmaller{N}}$. We partition $\g_{\mathsmaller{N-1}}$ as $g_{\mathsmaller{N-1}} = [\g_{\mathsmaller{N-2}}, g_{\mathsmaller{N}}]^{T}$. Now, Substituting the inner integral in (\ref{One/p1/Eqn/Ea32}), we have
%-----------------------------------
%------
\SeS
%-----------------------------------
%------------------------------------
\begin{EqA}{One/p1/Eqn/Ea33-1}{rcl}
%-----------------------------------
J &=&-2\INT{\re^{N-1}}{}{\phi(\zhat_{\mathsmaller{N-1}},\m_{\mathsmaller{N-1}},\R_{\mathsmaller{N-1}})\Phi\Bigg(\frac{\g_{\mathsmaller{N-1}}^{T}\zhat_{\mathsmaller{N-1}}-k_{\mathsmaller{N-1}}}{\sqrt{q_{\mathsmaller{N}}}}\Bigg)\exp(ju\zhats_1)\mathrm{d}\zhat_{\mathsmaller{N-1}}}
%-----------------------------------
\end{EqA}
%-----------------------------------
%-----------------------------------
%------------------------------------
\begin{EqA}{One/p1/Eqn/Ea33-2}{rcl}
%-----------------------------------
&=&-2\INT{\re^{N-1}}{}{\phi(\zhat_{\mathsmaller{N-2}},\m_{\mathsmaller{N-2}},\R_{\mathsmaller{N-2}}) \phi(\zhats_{\mathsmaller{N-1}},p_{\mathsmaller{N-1}},q_{\mathsmaller{N-1}})
\Phi\Bigg(\frac{g_{\mathsmaller{N-1}}\zhats_{\mathsmaller{N-1}}+\g_{\mathsmaller{N-2}}^{T}\zhat_{\mathsmaller{N-2}}-k_{\mathsmaller{N-1}}}{\sqrt{q_{\mathsmaller{N}}}}\Bigg)\exp(ju\zhats_1)\mathrm{d}\zhat_{\mathsmaller{N-1}}}\nonumber
%-----------------------------------
\end{EqA}
%-----------------------------------
%-----------------------------------
%------------------------------------
\begin{EqA}{One/p1/Eqn/Ea33-3}{rcl}
%-----------------------------------
&=&-2\INT{\re^{N-2}}{}{\Upsilon(\zhat_{\mathsmaller{N-2}})\Bigg(\INT{-\infty}{\infty}{\phi(\zhats_{\mathsmaller{N-1}}, p_{\mathsmaller{N-1}},q_{\mathsmaller{N-1}})\Phi\Bigg(\frac{g_{\mathsmaller {N-1}}\zhats_{N-1}+\g_{\mathsmaller{N-2}}^{T}\zhat_{\mathsmaller{N-2}}-k_{\mathsmaller{N-1}}}{\sqrt{q_{\mathsmaller{N}}}}\Bigg)\mathrm{d}\zhats_{\mathsmaller{N-1}}}\Bigg)\mathrm{d}\zhat_{\mathsmaller{N-2}}}\nonumber
%-----------------------------------
\end{EqA}
%-----------------------------------
\SeE
%-----------------------------------
where $\Upsilon(\zhat_{\mathsmaller{N-2}}) = \phi(\zhat_{\mathsmaller{N-2}},\m_{\mathsmaller{N-2}},\R_{\mathsmaller{N-2}})\exp(u\zhats_1)$, $p_{\mathsmaller{N-1}} = m_{\mathsmaller{N-1}}+\r^{T}_{\mathsmaller{N-1}}\R^{\dagger}_{\mathsmaller{N-2}}(\zhat_{\mathsmaller{N-2}}-\m_{\mathsmaller{N-2}})$ and $q_{\mathsmaller{N-1}} = r_{\mathsmaller{N-1N-1}}-\r^{T}_{\mathsmaller{N-1}}\R^{\dagger}_{\mathsmaller{N-2}}\r_{\mathsmaller{N-1}}$. To evaluate the inner integral, we make the transformation $\hat{t}_{\mathsmaller{N-1}}= \frac{\zhats_{\mathsmaller{N-1}} -p_{\mathsmaller{N-1}}}{\sqrt{q_{\mathsmaller{N-1}}}}$ to make it a standard integral, so we have
%-----------------------------------
%------------------------------------
\begin{EqA}{One/p1/Eqn/Ea34}{rcl}
%-----------------------------------
J&=&-2\INT{\re^{N-2}}{}{\Upsilon(\zhat_{\mathsmaller{N-2}})\Phi\Bigg(\frac{\alpha_{\mathsmaller{N-2}}}{\sqrt{1+\beta^2_{\mathsmaller{N-1}}}}\Bigg)\mathrm{d}\zhat_{\mathsmaller{N-2}}}
%-----------------------------------
\end{EqA}
%-----------------------------------
%-----------------------------------
where $\beta_{\mathsmaller{N-1}} = g_{\mathsmaller{N-1}}\sqrt{\frac{q_{\mathsmaller{N-1}}}{q_{\mathsmaller{N}}}}$ and $\alpha_{\mathsmaller{N-2}} = \frac{\g_{\mathsmaller{N-2}}^{T}\zhat_{\mathsmaller{N-2}}+g_{\mathsmaller{N-1}}p_{\mathsmaller{N-1}}
-k_{\mathsmaller{N-1}}}{\sqrt{q_{\mathsmaller{N}}}}$. Similarly evaluating the integral $N-3$ times, we have
%-----------------------------------
%------------------------------------
\begin{EqA}{One/p1/Eqn/Ea35}{rcl}
%-----------------------------------
J &=&-2\INT{-\infty}{\infty}{\phi(\zhats_{\mathsmaller{1}}, m_{\mathsmaller{1}},r_{\mathsmaller{11}})\Phi(\beta_{1}\zhats_{1}+\alpha_{1})\exp(ju\zhats_1)\mathrm{d}\zhats_{1}}
%-----------------------------------
\end{EqA}
%-----------------------------------
%-----------------------------------
which is equal to the Fourier transform of the function $\hat{f}(\zhats_1) = -2f(\zhats_1)\Phi(\beta_{1}\zhats_{1}+\alpha_{1})$, where $\alpha_1$ and $\beta_1$ are related to the entries of $\s$, $\m$ and $\R$. Hence,
$\mathcal{S}_{\u}[\cf_{\xhat}(\W\u)\star(\frac{-1}{\pi c_k})]$ is equal to the Fourier transform of $f(\zhats_1)+\hat{f}(\zhats_1) = f(\zhats_1)[1-2
\Phi(\beta_{1}\zhats_{1}+\alpha_{1})]$. So, we can observe from that when $\zhats_{k}$ is large and positive, then $\Phi(\beta\zhats_{k}+\alpha_{k})$ tends to zero and when $\zhats_{k}$ is large and negative, $\Phi(\beta\zhats_{k}+\alpha_{k})$ tends to one. So $(f(\zhats_k)+\hat{f}(\zhats_k))\approx f(\zhats_k)\sgn(\zhats_{k})$ for large $|\zhats_{k}|$, which justifies the use of this approximation in Corollaries \ref{One/p1/Thm/Thm6} and \ref{One/p1/Thm/Thm7}.
%-------------------------------------------------------
%\section*{Acknowledgements}
%---------------------------------------------
%---------------------------------------------
\bibliographystyle{plain}
\bibliography{references}

\begin{thebibliography}{10}

\bibitem{Mosesequivalence}
Christian~D. Austin, R.L. Moses, J.N. Ash, and E.~Ertin.
\newblock On the relation between sparse reconstruction and parameter
  estimation with model order selection.
\newblock {\em Selected Topics in Signal Processing, IEEE Journal of},
  4(3):560--570, 2010.

\bibitem{FL}
S.D. Babacan, R.~Molina, and A.K. Katsaggelos.
\newblock Bayesian compressive sensing using laplace priors.
\newblock {\em Image Processing, IEEE Transactions on}, 19(1):53--63, Jan 2010.

\bibitem{CS_Baraniuk}
R.G. Baraniuk, E.~Candes, R.~Nowak, and M.~Vetterli.
\newblock Compressive sampling [from the guest editors].
\newblock {\em Signal Processing Magazine, IEEE}, 25(2):12 --13, march 2008.

\bibitem{Ben-Tal-LMC}
Aharon Ben-Tal and Arkadiaei~Semenovich Nemirovskiaei.
\newblock {\em Lectures on Modern Convex Optimization: Analysis, Algorithms,
  and Engineering Applications}.
\newblock Society for Industrial and Applied Mathematics, Philadelphia, PA,
  USA, 2001.

\bibitem{cs_cv}
Petros Boufounos, Marco~F Duarte, and Richard~G Baraniuk.
\newblock Sparse signal reconstruction from noisy compressive measurements
  using cross validation.
\newblock In {\em Statistical Signal Processing, 2007. SSP'07. IEEE/SP 14th
  Workshop on}, pages 299--303. IEEE, 2007.

\bibitem{candes-rip}
E.~Candes.
\newblock {The restricted isometry property and its implications for compressed
  sensing}.
\newblock {\em Comptes Rendus Mathematique}, 346(9-10):589--592, May 2008.

\bibitem{chen-lasso}
Scott~Shaobing Chen, David~L. Donoho, and Michael~A. Saunders.
\newblock Atomic decomposition by basis pursuit.
\newblock {\em SIAM Rev.}, 43(1):129--159, January 2001.

\bibitem{cwold}
JuanAntonio Cuesta-Albertos, Ricardo Fraiman, and Thomas Ransford.
\newblock A sharp form of the cramér–wold theorem.
\newblock {\em Journal of Theoretical Probability}, 20(2):201--209, 2007.

\bibitem{donoho-cs}
D.L. Donoho.
\newblock Compressed sensing.
\newblock {\em Information Theory, IEEE Transactions on}, 52(4):1289--1306,
  April 2006.

\bibitem{Efronlars}
Bradley Efron, Trevor Hastie, Iain Johnstone, and Robert Tibshirani.
\newblock Least angle regression.
\newblock {\em Annals of Statistics}, 32:407--499, 2004.

\bibitem{sure-eldar}
Y.C. Eldar.
\newblock Generalized sure for exponential families: Applications to
  regularization.
\newblock {\em Signal Processing, IEEE Transactions on}, 57(2):471--481, Feb
  2009.

\bibitem{cvx}
Michael Grant and Stephen Boyd.
\newblock {CVX}: Matlab software for disciplined convex programming, version
  2.1.
\newblock \url{http://cvxr.com/cvx}, March 2014.

\bibitem{paul-critic}
Paul Kabaila.
\newblock The effect of model selection on confidence regions and prediction
  regions.
\newblock {\em Econometric Theory}, 11:537--549, 6 1995.

\bibitem{Kayesti}
Steven~M. Kay.
\newblock {\em Fundamentals of Statistical Signal Processing: Estimation
  Theory}.
\newblock Prentice-Hall, Inc., Upper Saddle River, NJ, USA, 1993.

\bibitem{knight2000}
Keith Knight and Wenjiang Fu.
\newblock Asymptotics for lasso-type estimators.
\newblock {\em Ann. Statist.}, 28(5):1356--1378, 10 2000.

\bibitem{Arraysp}
H.~Krim and M.~Viberg.
\newblock Two decades of array signal processing research: the parametric
  approach.
\newblock {\em Signal Processing Magazine, IEEE}, 13(4):67--94, Jul 1996.

\bibitem{potscher-critic}
Hannes Leeb and Benedikt~M. Pötscher.
\newblock Model selection and inference: Facts and fiction.
\newblock {\em Econometric Theory}, pages 21--59, 2 2005.

\bibitem{siglass}
Richard Lockhart, Jonathan Taylor, Ryan~J. Tibshirani, and Robert Tibshirani.
\newblock A significance test for the lasso.
\newblock {\em Ann. Statist.}, 42(2):413--468, 04 2014.

\bibitem{lopes-CS}
Miles~E. Lopes.
\newblock Estimating unknown sparsity in compressed sensing.
\newblock {\em CoRR}, abs/1204.4227, 2012.

\bibitem{Ng}
Ren Ng.
\newblock Fourier slice photography.
\newblock {\em ACM Trans. Graph.}, 24(3):735--744, July 2005.

\bibitem{complexlars}
A.~Panahi and M.~Viberg.
\newblock Fast candidate points selection in the lasso path.
\newblock {\em Signal Processing Letters, IEEE}, 19(2):79--82, Feb 2012.

\bibitem{potscher-finite}
Benedikt~M. P\"{o}tscher and Hannes Leeb.
\newblock On the distribution of penalized maximum likelihood estimators: The
  lasso, scad, and thresholding.
\newblock {\em J. Multivar. Anal.}, 100(9):2065--2082, October 2009.

\bibitem{tibs-lasso}
Robert Tibshirani.
\newblock Regression shrinkage and selection via the lasso.
\newblock {\em Journal of the Royal Statistical Society, Series B},
  58:267--288, 1994.

\bibitem{cvcs}
R.~Ward.
\newblock Compressed sensing with cross validation.
\newblock {\em Information Theory, IEEE Transactions on}, 55(12):5773--5782,
  Dec 2009.

\end{thebibliography}
%---------------------------------------------
\end{document}